\documentstyle[12pt]{amsart}

\numberwithin{equation}{section}

\begin{document}

%



\setlength{\itemsep}{0in}\newcommand{\lab}{\label}
\newcommand{\labeq}[1]{  \be \label{#1}  }
\newcommand{\labea}[1]{  \bea \label{#1}  }
\newcommand{\ben}{\begin{enumerate}}
\newcommand{\een}{\end{enumerate}}
\newcommand{\bm}{\boldmath}
\newcommand{\Bm}{\Boldmath}
\newcommand{\bea}{\begin{eqnarray}}
\newcommand{\ba}{\begin{array}}
\newcommand{\bean}{\begin{eqnarray*}}
\newcommand{\ea}{\end{array}}
\newcommand{\eea}{\end{eqnarray}}
\newcommand{\eean}{\end{eqnarray*}}
\newcommand{\beq}{\begin{equation}}
\newcommand{\eeq}{\end{equation}}
\newcommand{\bthm}{\begin{thm}}
\newcommand{\ethm}{\end{thm}}
\newcommand{\blem}{\begin{lem}}
\newcommand{\elem}{\end{lem}}
\newcommand{\bprop}{\begin{prop}}
\newcommand{\eprop}{\end{prop}}
\newcommand{\bcor}{\begin{cor}}
\newcommand{\ecor}{\end{cor}}
\newcommand{\bdfn}{\begin{dfn}}
\newcommand{\edfn}{\end{dfn}}
\newcommand{\brem}{\begin{rem}}
\newcommand{\erem}{\end{rem}}
\newcommand{\bex}{\begin{example}}
\newcommand{\eex}{\end{example}}
\newcommand{\lb}{\linebreak}
\newcommand{\nlb}{\nolinebreak}
\newcommand{\nl}{\newline}
\newcommand{\hs}{\hspace}
\newcommand{\vs}{\vspace}
\alph{enumii} \roman{enumiii}
\newtheorem{thm}{Theorem}[section]
\newtheorem{prop}[thm]{Proposition}
\newtheorem{lem}[thm]{Lemma}
\newtheorem{sublem}[thm]{Sublemma}
\newtheorem{cor}[thm]{Corollary}
\newtheorem{dfn}[thm]{Definition}
\newtheorem{rem}[thm]{Remark}
\newtheorem{example}[thm]{Example} 
\newtheorem{defn}{Definition} 

\def\endpf{\qed}
\def\ms{\medskip}
\def\N{{I\!\!N}}                \def\Z{Z\!\!\!\!Z}      \def\R{I\!\!R}
\def\C{{C\!\!\!\!I}}           \def\T{T\!\!\!\!T}      
\def\oc{\overline \C}
\def\Q{Q\!\!\!\!|}
\def\1{1\!\!1}

\def\and{\text{ and }}        \def\for{\text{ for }}
\def\tif{\text{ if }}         \def\then{\text{ then }}

\def\Cap{\text{Cap}}          \def\Con{\text{Con}}
\def\Comp{\text{Comp}}        \def\diam{\text{\rm {diam}}}
\def\dist{\text{{\rm dist}}}  \def\Dist{\text{{\rm Dist}}}
\def\Crit{\text{Crit}}
\def\Sing{\text{Sing}}        \def\conv{\text{{\rm conv}}}
\def\CIFS{\text{\rm CIFS}}    \def\PIFS{\text{\rm PIFS}}
\def\SIFS{\text{\rm SIFS}}
\def\HIFS{\text{\rm HIFS}}
\def\IR{\text{\rm IR}}        \def\RE{\text{\rm R}}
\def\CR{\text{\rm CR}}        \def\SR{\text{\rm SR}}
\def\CFR{\text{\rm CFR}}      \def\FSR{\text{\rm FSR}}

\def\RAH{\text{\rm RAH}}
\def\NPHH{\text{\rm NPHH}}

\def\CRIR{\text{\rm CRIR}}
\def\HSPT{\text{\rm HSPT}}

\def\Fin{{\cal F}in}
\def\F{{\Cal F}}
\def\h{{\text h}}
\def\hmu{\h_\mu}           \def\htop{{\text h_{\text{top}}}}

\def\H{\text{{\rm H}}}     \def\HD{\text{{\rm HD}}}   
\def\DD{\text{DD}}
\def\BD{\text{{\rm BD}}}         \def\PD{\text{PD}}
\def\re{\text{{\rm Re}}}    \def\im{\text{{\rm Im}}}
\def\Int{\text{{\rm Int}}} \def\ep{\text{e}}
\def\CD{\text{CD}}         \def\P{\text{{\rm P}}}     
\def\Id{\text{{\rm Id}}}   \def\Hol{\text{{\rm Hol}}}
\def\Hyp{\text{{\rm Hyp}}}
\def\A{\Cal A}             \def\Ba{\Cal B}       \def\Ca{\Cal C}
\def\Ha{{\cal H}}
\def\L{{\cal L}}             \def\M{\Cal M}        \def\Pa{\Cal P}
\def\U{\Cal U}             \def\V{\Cal V}  \def\Ra{{\mathcal R}}
\def\W{\Cal W}

\def\a{\alpha}                \def\b{\beta}             \def\d{\delta}
\def\De{\Delta}               \def\e{\varepsilon}          
\def\g{\gamma}                \def\Ga{\Gamma}                         
\def\l{\gamma}                \def\La{\Lambda}            
\def\Om{\Omega}               \def\om{\omega}
\def\Sg{\Sigma}               \def\sg{\sigma}
\def\Th{\Theta}               \def\th{\theta}           
\def\vth{\vartheta}
\def\Ka{\Kappa}               \def\ka{\kappa}

\def\bi{\bigcap}              \def\bu{\bigcup}
\def\({\bigl(}                \def\){\bigr)}
\def\lt{\left}                \def\rt{\right}
\def\bv{\bigvee}
\def\ld{\ldots}               \def\bd{\partial}         \def\^{\tilde}
\def\club{\hfill{$\clubsuit$}}\def\proot{\root p\of}

\def\tm{\widetilde{\mu}}
\def\tn{\widetilde{\nu}}
\def\es{\emptyset}            \def\sms{\setminus}
\def\sbt{\subset}             \def\spt{\supset}

\def\gek{\succeq}             \def\lek{\preceq}
\def\eqv{\Leftrightarrow}     \def\llr{\Longleftrightarrow}
\def\lr{\Longrightarrow}      \def\imp{\Rightarrow}
\def\comp{\asymp}
\def\upto{\nearrow}           \def\downto{\searrow}
\def\sp{\medskip}             \def\fr{\noindent}        
\def\nl{\newline}

\def\ov{\overline}            \def\un{\underline}

\def\ess{{\rm ess}}
\def\ni{\noindent}
\def\cl{\text{cl}}
\def\bt{{\bf t}}
\def\Bu{{\bf u}}
\def\tr{t}
\def\bo{{\bf 0}}
\def\nut{\nu_\lla^{\scriptscriptstyle 1/2}}
\def\arg{\text{arg}}
\def\Arg{\text{Arg}}
\def\re{\text{{\rm Re}}}
\def\gr{\nabla}
\def\endpf{{${\mathbin{\hbox{\vrule height 6.8pt depth0pt width6.8pt  
}}}$}}
\def\Fa{{\cal F}}
\def\Gal{{\cal G}}
\def\K{{\cal K}}
\def\1{1\!\!1}
\def\D{{I\!\!D}}

\title{Analytic Families of\\ Holomorphic Iterated Function Systems}
\author[\sc Mario ROY]{\sc Mario ROY}
\address{Mario Roy, Glendon College, York University, 
2275 Bayview Avenue, Toronto, Canada, M4N 3M6;
\newline {\tt mroy@@gl.yorku.ca
\newline Webpage: www.glendon.yorku.ca/mathematics/profstaff.html}}
\author[\sc Hiroki SUMI]{\sc Hiroki SUMI}
\address{Hiroki Sumi, Department of Mathematics, 
Graduate School of Science, Osaka University, 
1-1 Machikaneyama, Toyonaka, Osaka, 560-0043, Japan;
\newline {\tt sumi@@math.sci.osaka-u.ac.jp
\newline Webpage: www.math.sci.osaka-u.ac.jp/$\sim$sumi/}}
\author[\sc Mariusz URBA\'NSKI]{\sc Mariusz URBA\'NSKI}
\address{Mariusz Urba\'nski, Department of Mathematics,
University of North Texas, Denton, TX 76203-1430, USA;
\newline {\tt urbanski@@unt.edu
\newline Webpage: www.math.unt.edu/$\sim$urbanski}}
%
\date{August 16, 2008. Published in Nonlinearity 21 (2008) 2255-2279.}
\thanks{Research of the first author was supported by NSERC 
(Natural Sciences and Engineering Research Council of Canada).
The research of the third author was supported in part by the
NSF Grant DMS 0700831.}

\begin{abstract}
This paper deals with analytic families of holomorphic iterated function 
systems. Using real analyticity of the pressure function (which we prove), we 
establish a classification theorem for analytic families of holomorphic iterated function 
systems which depend continuously on a parameter when the space of holomorphic iterated function 
systems is endowed with the $\lambda$-topology. 
This classification theorem allows us to generalize some geometric results from ~\cite{RSU} and 
gives us a better and clearer understanding of the global structure of the space of conformal IFSs.
\end{abstract}

\maketitle
\section{\bf Introduction}

Iterated function systems (abbreviated to IFSs henceforth) arise in many 
natural contexts. They are often used to encode and generate fractal 
images, such as landscapes and skyscapes, in computer games. They each 
generate, via a recursive procedure, a unique fractal set called attractor or limit set. 
IFSs also play an important role in the theory of dynamical systems. 
By dynamical system we mean a continuous map $T$ from a metric space $X$ to itself, 
where, given $x \in X$, one aims at describing the eventual behavior of the 
sequence of iterates $(T^{n}x)_{n=0}^\infty$. IFSs are in fact a generalization 
of the process of looking at the backward trajectories of dynamical systems. 
A repeller in a dynamical system sometimes coincides with the attractor of an 
associated IFS. For example, the middle-third Cantor set attracts all the 
backward orbits of the tent map: $T(x)=3x$ if $x \leq 1/2$ and $T(x)=3(1-x)$ 
if $x \geq 1/2$. This follows from the fact that the inverse branches of $T$ 
are precisely the generators $\varphi_{1}$ and $\varphi_{2}$ of an IFS. 
As dynamical systems often model physical processes, the study of IFSs 
frequently turns out to be instrumental in describing 
real systems. 

The systematic development of modern theory of iterated function systems 
 began 
with the works of Hutchinson~\cite{Hu}, Falconer~\cite{Fal}, Barnsley~\cite{Barn}, 
Bandt and Graf~\cite{BG}, and Schief~\cite{Sch}, just to name a few. 
Their pioneering works in the 1980's and 1990's on finite IFSs concerned 
systems consisting of similarities. In particular, the theory of finite IFSs 
was used in the study of the complex dynamics of rational functions. 

In the middle of the 1990's, the need to investigate finite, 
and even infinite, IFSs that comprise more general conformal 
maps arose. The foundations of that theory were laid out by Mauldin 
and Urba\'nski in \cite{MU}. Many new applications were also discovered. In 1999,
Mauldin and Urba\'nski~\cite{mutr} applied the theory of infinite conformal IFSs to continued 
fractions with restricted entries. Few years later, 
Kotus and Urba\'nski~\cite{KU} applied that theory to obtain a lower estimate 
for the Hausdorff dimension of the Julia set of elliptic functions, while Urba\'nski and 
Zdunik~\cite{UZ} obtained similar results in the complex dynamics of the exponential function.
They further showed that the framework of conformal IFSs is the right one to study the 
harmonic measure of various Cantor sets. Applications of conformal IFSs to number theory, 
more precisely to the theory of continued fractions and 
Diophantine approximations, were also developed by Urba\'nski in~\cite{U1}, 
\cite{U2}, \cite{U3}, and \cite{U4}. Finally, Stratmann and Urba\'nski~\cite{SU0} 
established extremality in the sense of Kleinbock, Lindenstrauss and Barak Weiss, 
of conformal measures for convex co-compact Kleinian groups.

A further generalization of the theory was achieved by Mauldin and Urba\'nski~\cite{gdms} in 2003. This 
generalization relies on graphs. Indeed, every IFS can be thought of as a graph with a single vertex 
and a countable set of self-loops. To the unique vertex is attached a 
space $X\sbt\R^d$. The self-loops represent the generators of the IFS, each generator being 
a contracting self-map of $X$. The limit set generated by all infinite paths on 
this graph is also attached to the vertex. 
Mauldin and Urba\'nski extended the theory of IFSs 
to graph directed Markov systems (GDMSs). The graphs associated with these
systems generally have more than one, though finitely many, vertices. Moreover, they can have a countably 
infinite set of edges between any two of their vertices. It turns out that a unique fractal set can be 
associated with each vertex, and the limit set is the disjoint union of those sets. In 2007,
Stratmann and Urba\'nski~\cite{SU} went one step further by laying the foundations to the theory of pseudo-Markov systems. In contradistinction with GDMSs, the underlying graph of these systems can have infinitely many vertices. 
They applied their work to infinitely generated Schottky groups. 

In all those works, the metric structure of the limit set $J$ is usually described by its Hausdorff, packing and/or box-counting dimension(s). The Hausdorff dimension is particularly interesting as it is characterized by the pressure function. The pressure function plays a central role in thermodynamic formalism, which arose from statistical
physics. For finite systems which satisfy a certain separation property (the famous open set condition (OSC)), 
the Hausdorff dimension of the limit set is simply the zero of the pressure 
function $P(t)$, that is, the unique $t_0>0$ such that $P(t_0)=0$. This latter equation is sometimes 
called Moran-Bowen formula. For infinite systems satisfying the OSC and a bounded distortion property, Mauldin and Urba\'nski~\cite{MU} showed that a variant of Moran-Bowen formula holds: the Hausdorff dimension of the limit set is  the infimum of all $t\ge 0$ for which $P(t)\le 0$. They also showed that a $t$-conformal measure exists 
if and only if $P(t)=0$, and there is at most one such $t$.

Recently, interest in families of IFSs and GDMSs has emerged 
(see \cite{As},~\cite{BR},~\cite{RU},~\cite{RSU} and~\cite{RURA}, among others). 
Baribeau and Roy~\cite{BR} showed that the Hausdorff dimension of the 
limit set is a continuous, subharmonic function in families of (finite or infinite) 
IFSs which consist of complex similarities whose coefficients vary holomorphically with a 
parameter. Using thermodynamic formalism, Roy and Urba\'nski~\cite{RU} later 
proved that this result extends to families of conformal IFSs. 
Moreover, they produced diverse families of infinite IFSs in which the Hausdorff dimension 
of the limit set does not vary real-analytically: there are phase transitions during 
which the real-analyticity sometimes breaks 
down. This is in sharp contrast to what happens in the finite case.
In a subsequent paper, Roy and Urba\'nski~\cite{RURA} studied 
the behavior of the Hausdorff dimension function 
of the limit sets of conformal GDMSs that live 
in higher dimensional Euclidean spaces $\R^d$, $d\ge 3$. They proved in particular that this function 
is real-analytic in the subspace of all strongly regular systems that have a finitely irreducible 
incidence matrix.  

Roy and Urba\'nski also considered a slightly different approach 
in~\cite{RU}. They fixed an underlying space $X\sbt\R^d$ and a
countable index set $I$, and investigated the space $\mbox{CIFS}(X,I)$ of 
all conformal IFSs with alphabet $I$ that live on $X$. 
When $I$ is finite and $\CIFS(X,I)$ is endowed with the 
natural metric of pointwise convergence (pointwise meaning 
that corresponding generators are compared to one another in $C^1(X)$), 
the topological pressure and the Hausdorff dimension functions 
are continuous functions of the underlying CIFSs. 
However, they discovered that these functions 
are generally not 
continuous when $I$ is infinite (and $\CIFS(X,I)$ is equipped with a 
metric of pointwise convergence). 
They thereafter introduced a finer but still natural topology 
called the $\lambda$-topology, 
with respect to which the topological pressure and the Hausdorff 
dimension functions are both continuous (see~\cite{RU}). 

More recently, the authors of the present article have studied 
more thoroughly the pointwise and $\lambda$- topologies (see~\cite{RSU}). We have shown that
the $\lambda$-topology is normal but not metrizable, for it is not first-countable. 
We studied in detail how various characteristics of an IFS behave with respect to the $\lambda$-topology.
 
The aim of the current paper is to investigate some families in a subspace of $\CIFS(X,I)$ 
when $d=2$ (or, equivalently, when $X\sbt\C$) and when $\CIFS(X,I)$ is endowed 
with the $\lambda$-topology. The abovementioned subspace, denoted by $\HIFS(X,I)$, 
is the space of all holomorphic CIFSs. We shall study analytic families 
of HIFSs whose generators depend analytically on a complex parameter. 

After some preliminaries on single CIFSs in section~\ref{secprelifs} and on 
families of CIFSs in section~\ref{fams}, we study in section~\ref{pf} the properties of the 
pressure function seen as a function of two variables, that is, as a function of not only 
the usual parameter $t$ but also of the underlying IFS $\Phi$. In section~\ref{ct}, we enunciate 
a classification theorem of analytic families which depend continuously on a parameter when 
$\HIFS(X,I)$ is endowed with the $\lambda$-topology. Finally, in section~\ref{cct}, we use 
our classification theorem to generalize some results that we obtained in~\cite{RSU} and hence 
to gain a better understanding of $\CIFS(X,I)$ when this latter is equipped with the 
$\lambda$-topology. 

 Furthermore, we give some examples of analytic families 
 in HIFS$(X,I)$ (Example~\ref{example}, Proposition~\ref{propo}, Theorem~\ref{t:hsptdense}, 
 Proposition~\ref{p:everytype}).  
\

\section{Preliminaries on Iterated Function Systems}\label{secprelifs}
Let us first describe the setting of conformal iterated function systems 
introduced in \cite{MU}. Let $I$ be a countable (finite or infinite) index 
set (often called alphabet) with at least two elements, and let 
$\Phi=\{\varphi_i:X\to X \, | \, i\in I\}$ be a collection 
of injective contractions of a compact metric space $(X,d_X)$ 
(sometimes coined seed space) for which there exists a constant $0<s<1$ such
that $d_X(\varphi_i(y),\varphi_i(x)) \leq s \, d_X(y,x)$ for every $x,y\in X$
and for every $i\in I$. Any such
collection $\Phi$ is called an iterated function
system (abbr. IFS). We define the limit set $J_\Phi$ of this system as the 
image of the coding space $I^\infty$ under a coding map $\pi_\Phi$,
which is defined as follows.
Let $I^n$ denote the space of words of length $n$ with letters in $I$, let $I^*:=\bigcup_{n\in\N}I^n$ 
be the space of finite words, and $I^\infty$ the space of
one-sided infinite words with letters in $I$. 
Given $\omega, \tau \in I^{\infty }$, we define $\omega \wedge \tau \in I^* \cup I^{\infty }$ 
to be the longest initial block common to both $\omega $ and $\tau .$  
For every $\om \in I^*\cup I^\infty$, we write $|\om|$ for the length 
of $\om$, that is, the unique $n \in \N\cup \{\infty\}$ such that $\om\in I^n$.
For $\om\in I^n$ with $n\in\N$, we set 
$\varphi_\om:=\varphi_{\om_1}\circ\varphi_{\om_2}\circ\cdots\circ \varphi_{\om_n}$. 
If $\om\in I^*\cup I^\infty$ 
and $n\in\N$ does not exceed the length 
of $\om$, we denote by $\om|_n$ the word $\om_1\om_2\ldots\om_n$. Since, 
given $\om\in I^\infty$, the diameters of the compact sets $\varphi_{\om|_n}(X)$,
$n\in\N$, converge to zero and since these sets form a decreasing family, their 
intersection
$$
\bigcap_{n=1}^\infty\varphi_{\om|_n}(X)
$$
is a singleton, and we denote its element by $\pi_\Phi(\om)$. This 
defines the coding map $\pi_\Phi:I^\infty\to X$. Clearly, $\pi_\Phi$ is a continuous 
function when $I^\infty$ is equipped with the topology generated by the 
cylinders $[i]_{n}=\{\om\in I^\infty \, : \, \om_n=i\}$, $i\in I$, $n\in\N$. 
The main object of our interest is the limit set
$$
J_\Phi=\pi_\Phi(I^\infty)=\bigcup_{\om\in I^\infty}\bigcap_{n=1}^\infty\varphi_{\om|_n}(X).
$$
Observe that $J_\Phi$ satisfies the natural invariance property, $J_\Phi=
\bigcup_{i\in I}\varphi_i(J_\Phi)$. Note also that if $I$ is finite, then $J_\Phi$ is
compact, which is usually not the case when $I$ is infinite.

\ms\ni An IFS $\Phi=\{\varphi_i\}_{i\in I}$ 
is said to satisfy the Open Set Condition (OSC) if there exists 
a nonempty open set $U \sbt X$ such that $\varphi_i(U)\sbt U$ for every $i\in I$ 
and $\varphi_i(U)\cap\varphi_j(U)=\emptyset$ for every pair $i,j\in I$, $i\neq j$.

\ms\ni An IFS $\Phi$ is called conformal (and thereafter a CIFS) if $X$ is 
connected, $X=\overline{\Int_{\R^d}(X)}$ for some $d\in\N$, and if the following 
conditions are satisfied:
\begin{itemize}
\item[(i)] $\Phi$ satisfies the OSC with $U=\Int_{\R^d}(X)$;
\item[(ii)] There exists an open connected set $V$, with 
$X\sbt V\sbt
\R^d$, such that all the generators $\varphi_i$, $i\in I$, extend to $C^1$ 
conformal diffeomorphisms of $V$ into $V$;
\item[(iii)] There exist $\g,l>0$ such that for every $x\in X$ 
there is an open cone $\Con(x,\g,l)\sbt \Int(X)$ with vertex $x$, central 
angle of Lebesgue measure $\g$, and altitude $l$;
\item[(iv)] There are two constants $L\ge 1$ and $\a>0$ such that
$$
\bigl||\varphi_i'(y)|-|\varphi_i'(x)|\bigr|\le L\|(\varphi_i')^{-1}\|_V^{-1}\|y-x\|^\a
$$
for all $x,y\in V$ and all $i\in I$, where $\|\cdot\|_V$ is the supremum norm taken over $V$.
\end{itemize}

\

\brem\lab{rsbdp}
It has been proved in Proposition~4.2.1 of \cite{gdms} that if $d\ge 2$,
then condition $(iv)$ is satisfied with $\a=1$. This condition is also 
fulfilled if $d=1$, the alphabet $I$ is finite and all the maps
$\varphi_i$ are of class $C^{1+\e}$.
\erem

\

\ni The following useful fact has been also proved in~Lemma~4.2.2 of~\cite{gdms}.

\

\blem\lab{lsbdp}
For all $\om\in I^*$ and all $x,y\in V$ we have that
$$
\bigl|\log|\varphi_\om'(y)|-\log|\varphi_\om'(x)|\bigr|
\le L(1-s)^{-1}\|y-x\|^\a.
$$
\elem

\

\ni As an immediate consequence of this lemma we get the following.

\

\begin{itemize}
\item[(iv')] Bounded Distortion Property (BDP): There exists a constant $K\geq
1$ such that
$$
|\varphi_\om'(y)|\leq K|\varphi_\om'(x)|
$$
for every $x,y\in V$ and every $\om\in I^*$, where $|\varphi_\om'(x)|$ 
denotes the norm of the derivative.
\end{itemize}


\

\fr As demonstrated in \cite{MU}, infinite CIFSs, unlike finite ones,
split naturally into two main classes: irregular and regular systems. This
dichotomy can be expressed in terms of the absence or existence of a zero for the 
topological pressure function. Recall that the topological pressure  
$\P_\Phi(t)$, $t\ge 0$, is 
defined as follows. For every $n\in\N$, set
$$
\P_\Phi^{(n)}(t)=\sum_{\om\in I^n}\|\varphi_\om'\|^t,
$$
where $\|\cdot\|:=\|\cdot\|_X$ is the supremum norm over $X$.
Then
$$
\P_\Phi(t)=\lim_{n\to\infty}{1\over n}\log \P_\Phi^{(n)}(t)=\inf_{n\in\N}{1\over n}\log \P_\Phi^{(n)}(t).
$$
Recall also that the shift map
$\sg:I^*\cup I^\infty\to I^*\cup I^\infty$ is defined for each $\om\in I^*\cup I^\infty$ as 
$$
\sg\(\{\om_n\}_{n=1}^{|\om|}\)=\{\om_{n+1}\}_{n=1}^{|\om|-1}.
$$
Now, if the function $\zeta_\Phi:I^\infty\to\R$ is given by
the formula
$$
\zeta_\Phi(\om)=\log|\varphi_{\om_1}'(\pi(\sg(\om)))|,
$$
then $\P_\Phi(t)=\P(t\zeta_\Phi)$, where $\P(t\zeta_\Phi)$ is the classical 
topological pressure of the function $t\zeta_\Phi$ when $I$ is finite 
(so the space $I^\infty$ is compact), and is understood 
in the sense of \cite{HMU} and \cite{gdms} when $I$ is infinite.
The finiteness parameter $\th_\Phi$ of a system 
is defined as $\inf\{t\ge 0:\P_\Phi^{(1)}(t)<\infty\}=
\inf\{t\ge 0:\P_\Phi(t)<\infty\}$. Let also
\[ \mbox{Fin}(\Phi)=\bigl\{t\ge 0:P_\Phi(t)<\infty\bigr\}=\bigl\{t\ge 0:P_\Phi^{(1)}(t)<\infty\bigr\}. \] 
In \cite{MU}, it was shown that 
the topological pressure function $\P_\Phi$ is non-increasing on $[0,\infty)$, 
(strictly) decreasing, continuous and convex on 
$[\th_\Phi,\infty)$, and $\P_\Phi(d) \leq 0$. Of course, $\P_\Phi(0)=\infty$ if and 
only if $I$ is infinite.
The following characterization of the Hausdorff dimension 
$h_\Phi$ of the limit 
set $J_\Phi$ was stated as Theorem~3.15 in~\cite{MU}. 
For every $F\sbt I$, we write $\Phi|_F$ for the 
subsystem $\{\varphi_i\}_{i\in F}$ of $\Phi$.

\

\bthm\lab{3.15mul}
$$
h_\Phi = \sup\{h_{\Phi|_F}:F\sbt I  \text{ is finite }\}
= \inf\{t\ge 0:\P_\Phi(t)\le 0\}.
$$
In particular, if $\P_\Phi(t)=0$, then $t=h_\Phi$.
\ethm

\

\fr Subsequently, a system $\Phi$ was called regular provided there is some $t\ge 0$ 
such that $\P_\Phi(t)=0$. It follows from the strict decrease of $P_\Phi$ on $[\th_\Phi,\infty)$ 
that such a $t$ is unique, if it exists. If no such $t$ exists, the system is called irregular. 

\sp\fr Regular systems can also be naturally divided into subclasses. Following~\cite{MU} still, 
a system $\Phi$ is said to be strongly regular if 
$0<\P_\Phi(t)<\infty$ for some $t\ge 0$. As an immediate application of 
Theorem~\ref{3.15mul}, a system $\Phi$ is strongly regular 
if and only if $h_\Phi>\th_\Phi$. A system $\Phi$ is called  
cofinitely regular provided that every nonempty 
cofinite subsystem $\Phi' = \{\varphi_i\}_{i\in I'}$ (with $I'$ any cofinite 
subset of $I$) is regular.
A finite system is clearly cofinitely regular, and it 
was shown in \cite{MU} that an infinite system is cofinitely regular exactly when  
the pressure is infinite at the finiteness parameter, that is, $\P_\Phi(\th_\Phi)=\infty$.
Note that every cofinitely regular system is strongly regular, and every 
strongly regular system is regular. Finally, recall from~\cite{RU} that critically regular systems 
are those regular systems for which $P_\Phi(\th_\Phi)=0$.

\section{Preliminaries on Families of CIFSs}\label{fams}

When dealing with families of CIFSs, we denote the set of all conformal 
iterated function systems with phase space $X$ and alphabet $I$ by 
$\CIFS(X,I)$. Obviously, $\CIFS(X,I)$ can be endowed with different topologies. 

When $I$ is finite, $\CIFS(X,I)$ is naturally endowed with the metric of 
pointwise convergence. This metric asserts that the distance between 
CIFSs $\Phi=\{\varphi_i\}_{i\in I}$ and $\Psi=\{\psi_i\}_{i\in I}$ is 
$$
\rho(\Phi,\Psi)
=\sum_{i\in I}\(\|\varphi_i-\psi_i\| + \|\varphi_i'-\psi_i'\|\).
$$
It was proved in~\cite{RU} that, given $t\ge 0$, the pressure function
$P(t):\CIFS(X,I)\to\R$, $\Phi\mapsto P_\Phi(t)$ and the Hausdorff dimension 
function $h:\CIFS(X,I)\to[0,d]$, $\Phi\mapsto h_\Phi$ are continuous 
(see Lemma~4.2 and Theorem~4.3 in \cite{RU}, respectively). 

When $I$ is infinite, the situation is more intricate. First, recall that 
we may assume that $I=\N$ without loss of generality. 
Henceforth, we accordingly abbreviate $\CIFS(X,\N)$ to $\CIFS(X)$. 
The set $\CIFS(X)$ 
can easily be equipped with a metric of pointwise convergence. 
In~\cite{RU}, such a metric was introduced by defining 
the distance between two CIFSs $\Phi$ and $\Psi$ as  
$$
\rho_\infty(\Phi,\Psi)
=\sum_{i=1}^\infty 
2^{-i}\min\bigl\{1,\|\varphi_i-\psi_i\|+\|\varphi_i'-\psi_i'\|\bigr\}.
$$
In these definitions, $\|\cdot\|:=\|\cdot\|_X$ stands for the supremum norm over $X$. 

Roy and Urba\'nski observed that the pressure and Hausdorff dimension 
functions are generally not continuous when $\CIFS(X)$ is endowed with
the metric $\rho_\infty$ 
(see the example following Lemma~5.3 in~\cite{RU}). This triggered the 
introduction of a topology called the $\lambda$-topology (see section~5 in~\cite{RU}). 
In that topology, a sequence $\{\Phi^n\}$ converges to 
$\Phi$ provided that $\{\Phi^n\}$ converges to $\Phi$ 
in the pointwise topology and that there 
exist constants $C>0$ and $N\in\N$ such that
\begin{equation}\label{1061804}
\big|\log\|\varphi_i'\|-\log\|(\varphi_i^n)'\|\big|\le C
\end{equation}
for all $i\in\N$ and all $n\ge N$. A set $F\sbt\CIFS(X)$ is declared to 
be closed if the $\lambda$-limit of every $\lambda$-converging sequence of CIFSs
in $F$ belongs to $F$. Several topological properties of $\CIFS(X)$ 
were given in~\cite{RU} and~\cite{RSU}. Among others, let us mention that 
this topology is not metrizable, for it does not even satisfy the first 
axiom of countability (see Proposition~5.7 in~\cite{RSU}). Nonetheless, 
this topology proves to be useful, for it is easy to determine whether 
a sequence converges or not in that topology and, according to Roy and 
Urba\'nski, the Hausdorff dimension function is 
then continuous everywhere on $\CIFS(X)$ (see Theorem~5.10 in~\cite{RU}). 
In fact, the combination of Theorem~5.7 in~\cite{RU} with Lemma~5.22 
in~\cite{RSU} shows that the pressure function is continuous wherever it 
possibly can be.

Whichever topology we choose to endow $\CIFS(X)$ with, there are some subspaces
of particular interest. Following the notation in~\cite{RSU}, let $\SIFS(X)$ 
represent the subset of $\CIFS(X)$ comprising all similarity iterated function 
systems, that is, iterated function systems whose generators are similarities. 
Let also $\IR(X)\sbt\CIFS(X)$ be the subset of irregular systems, while 
$\RE(X)\sbt\CIFS(X)$ will represent the subset of regular systems. 
Denote further by $\CR(X)\sbt\RE(X)$ the subset of critically regular systems, by 
$\SR(X)\sbt\RE(X)$ the subset of strongly regular systems, and by 
$\CFR(X)\sbt\SR(X)$ the subset of cofinitely regular systems. Finally, we will
denote by $\FSR(X)$ the set $\SR(X)\sms\CFR(X)$.

In this paper, we will concentrate on the case $d=2$ or, 
equivalently, $X\sbt\C$. In this case, 
a natural subset of $\CIFS(X,I)$ is: 

\bdfn
The set $\HIFS(X,I)$ consists of those systems $\Phi=\{\varphi_i\}_{i\in I}
\in\CIFS(X,I)$ which admit an open   
connected neighborhood $V$ of $X\sbt\C$ such that
for each $i\in I$, the map $\varphi _{i}$ 
extends to a holomorphic diffeomorphism of $V$ into $V.$ 
\edfn

\

More specifically, we will be interested in analytic families of HIFSs. 

\bdfn\label{df}
A family $\{\Phi^\g\}_{\g\in\Ga}
=\{\{\varphi_i^\g\}_{i\in I}\}_{\g\in\Ga}$ 
in $\HIFS(X,I)$ is called analytic if
\begin{itemize} 
\item $\Ga$ is a connected, finite-dimensional complex manifold; and
\item for each 
$\g_0\in\Ga $,  
there exists a neighborhood $\Ga_0$ of $\g_0$ and 
a bounded connected open neighborhood $V_0$ of $X\subset\C$ 
such that for each $i\in I$ and each $\g \in \Ga $, 
$\varphi _{i}^{\g }$ extends to a holomorphic diffeomorphism 
of $V_{0}$ into $V_{0}$, and for each $i\in I$, 
the map 
\[(\g,x)\mapsto\varphi_i^\g(x), \hspace{0.5cm} 
(\g,x)\in\Ga_0\times V_0, \]
is holomorphic.  
\end{itemize}
\edfn

\section{The Pressure as a Function of Two Variables}\label{pf}

In this section, we make some observations about the pressure function seen
as a function of two variables, that is, as a function of not only the parameter 
$t$ but also of the underlying CIFS $\Phi$. These observations mostly follow 
from earlier results presented in~\cite{RU} and~\cite{RSU}. 

When the alphabet $I$ is finite, the pressure behaves well. 

\begin{lem}\label{l1042903}
If $I$ is a finite alphabet, then the 
pressure function $(\Phi,t)\mapsto P(\Phi,t)=P(t\zeta_\Phi)$,
$\Phi\in\CIFS(X,I)$, $t\geq 0$, is continuous. 
\end{lem}

{\sl Proof.}  In the proof of Lemma~4.2 in~\cite{RU}, it was shown that 
the function $\Phi\mapsto \zeta_\Phi\in C(I^\infty)$ is continuous. 
Since the pressure function $P:C(I^\infty)\to\R$ is
Lipschitz continuous with Lipschitz constant $1$, we deduce that
\begin{eqnarray*}
\bigl|P(\Phi,t)-P(\Phi_0,t_0)\bigr|
&=&\bigl|P(t\zeta_\Phi)-P(t_0\zeta_{\Phi_0})\bigr| \\
&\leq&\|t\zeta_\Phi-t_0\zeta_{\Phi_0}\| \\
&\leq&\|t\zeta_\Phi-t\zeta_{\Phi_0}\|+\|t\zeta_{\Phi_0}-t_0\zeta_{\Phi_0}\| \\
&=&|t|\cdot\|\zeta_\Phi-\zeta_{\Phi_0}\|+|t-t_0|\cdot\|\zeta_{\Phi_0}\|,
\end{eqnarray*}
for every $\Phi,\Phi_0\in\CIFS(X,I)$ and $t,t_0\geq 0$.
Thus, the pressure function $(\Phi,t)\mapsto P(\Phi,t)$ is continuous.
\endpf

\

For an infinite alphabet, we obtain immediately the following. Recall that for any $\Phi\in\CIFS(X)$, 
\[ \mbox{Fin}(\Phi)=\bigl\{t\ge 0:P(\Phi,t)<\infty\bigr\}. \]  

\begin{lem}\label{l1042903c}
The pressure function $(\Phi,t)\mapsto P(\Phi,t)$,
$\Phi\in\CIFS(X)$, $t\geq 0$, is lower semi-continuous when $\CIFS(X)$ 
is endowed with the metric $\rho_\infty$ of pointwise convergence. 
Moreover, it is continuous on the set $\cup_{\Phi\in\CIFS(X)}
\{\Phi\}\times Fin(\Phi)^c$, and discontinuous on $\cup_{\Phi\in\CIFS(X)}
\{\Phi\}\times (Fin(\Phi)\cap[0,d))$. In fact, in this latter case,
given any $\Phi\in\CIFS(X)$ and $t\in Fin(\Phi)\cap[0,d)$, 
the function $\Psi\mapsto P(\Psi,t)$ is discontinuous at $\Phi$. 
\end{lem}

{\sl Proof.}  Using Theorem~2.1.5 in~\cite{gdms} as well as 
Lemma~\ref{l1042903}, for any finite subset $F$ of $\N$ we have that
\[ \liminf_{(\Phi,t)\to(\Phi_0,t_0)} P(\Phi,t) 
\ge \liminf_{(\Phi,t)\to(\Phi_0,t_0)} P_F(\Phi,t) 
= \lim_{(\Phi,t)\to(\Phi_0,t_0)} P_F(\Phi,t)
= P_F(\Phi_0,t_0). \] 
Since this is true for every finite subset $F$ of $\N$, we deduce from 
Theorem~2.1.5 in~\cite{gdms} 
that
\[ \liminf_{(\Phi,t)\to(\Phi_0,t_0)} P(\Phi,t)\geq
\sup_{F\subset\N, F {\tiny \mbox{finite}}} P_F(\Phi_0,t_0) 
= P(\Phi_0,t_0). \] 

Observe also that $P(\Phi,t)=\infty$ for any $t\notin Fin(\Phi)$ and 
every $\Phi\in\CIFS(X)$. Thus, the pressure function is upper semi-continuous
at any point in the set $\cup_{\Phi\in\CIFS(X)}
\{\Phi\}\times Fin(\Phi)^c$.

Finally, given $\Phi\in\CIFS(X)$ and $t\in Fin(\Phi)\cap[0,d)$, 
the function $\Psi\mapsto P(\Psi,t)$ 
is discontinuous at $\Phi$ according to the last paragraph in the proof of
Lemma~4.4 in~\cite{RSU}.
\endpf

\

When $\CIFS(X)$ is endowed with the $\lambda$-topology, 
we have the following result.

\begin{lem}\label{l1042903cl}
The pressure function $(\Phi,t)\mapsto P(\Phi,t)$
is continuous on the set $\cup_{\Phi\in\CIFS(X)}
\{\Phi\}\times Fin(\Phi)^c$ and sequentially continuous on
$\cup_{\Phi\in\CIFS(X)}\{\Phi\}\times Fin(\Phi)$
when $\CIFS(X)$ 
is endowed with the $\lambda$-topology. 
\end{lem}

{\sl Proof.}  The continuity of the pressure on the set $\cup_{\Phi\in\CIFS(X)}
\{\Phi\}\times Fin(\Phi)^c$ follows from Lemma~\ref{l1042903c}. The sequential 
continuity on $\cup_{\Phi\in\CIFS(X)}\{\Phi\}\times Fin(\Phi)$ follows from 
Lemma~5.22 in~\cite{RSU} and the fact that if $\{\Phi^n\}$ converges to 
$\Phi$ in the $\lambda$-topology, then $Fin(\Phi^n)=Fin(\Phi)$ for all $n$ 
large enough, as observed before Theorem~5.20 in~\cite{RSU}. 
\endpf 

\

We now investigate the properties of the pressure function 
$(\g,t)\mapsto P(\g,t)=P(\Phi^\g,t)$ 
within analytic families in $\HIFS(X,I)$. 

We first prove that every analytic family in $\HIFS(X,I)$ constitutes a continuous family 
in $\CIFS(X,I)$ when this latter is equipped with the pointwise topology.

\

\begin{thm}\label{l1042903t}
An analytic family $\{\Phi^\g\}_{\g\in\Ga}$ regarded 
as the function 
\[ \begin{array}{ccccc} \Phi & : & \Ga & \rightarrow &\CIFS(X,I) \\
                             &   & \g & \mapsto     & \Phi^\g
   \end{array} 
\] 
is continuous when $\CIFS(X,I)$ is endowed with the metric of pointwise 
convergence ($\rho$ or $\rho_\infty$ depending on the case). 
\end{thm}

{\sl Proof.} 
The proof in the infinite case reduces to the finite case. Indeed, if 
$\Phi$ were discontinuous at some $\g_0\in\Ga$, there would 
exist $\e>0$ and a sequence $\{\g_n\}$ such that 
$\g_n\rightarrow\g_0$ but 
$\rho_\infty(\Phi^{\g_n},\Phi^{\g_0})\geq\e$ 
for every $n\in\N$. Choose $J\in\N$ such that 
$\sum_{j>J}2^{-j}<\e/2$. Then 
\[ \sum_{j\leq J}2^{-j}\min\bigl\{1,\|\varphi_j^{\g_n}-\varphi_j^{\g_0}\|+
\|(\varphi_j^{\g_n})'-(\varphi_j^{\g_0})'\|\bigr\}>\e/2 \]
since, otherwise, we would have $\rho_\infty(\Phi^{\g_n},\Phi^{\g_0})<\e$. 
Replacing the sequence $\{\g_n\}$ by one of its subsequences if necessary, 
it follows from this that there is an $i\leq J$ such that 
\[ \|\varphi_i^{\g_n}-\varphi_i^{\g_0}\|+
\|(\varphi_i^{\g_n})'-(\varphi_i^{\g_0})'\|>\frac{\e}{2J} \]
for every $n\in\N$. Consequently, each subsequence of $\{\varphi_i^{\g_n}\}_{n\in\N}$ either does not converge uniformly on $X$ to $\varphi_i^{\g_0}$ or the corresponding subsequence of derivatives does not converge uniformly on $X$ to $(\varphi_i^{\g_0})'$.

However, given the neighborhood $V_{0}$ from the definition of an analytic family, the generators $\{\varphi_j^\g\}_{j\in\N, \g\in\Ga}$ form a family of holomorphic maps uniformly bounded on $V_{0}$, say by a constant $M$, due to the forward invariance of $V_{0}$. By Cauchy's Integral Formula, their derivatives are thereby uniformly bounded on any compact subset $K$ of $V_{0}$ by a constant that depends on $M$, $K$ and $V_{0}$ only (but neither on $j$ nor on $\g$). Therefore the generators constitute an equicontinuous family on $V_{0}$ and thereby form a normal family by Arzel\`a-Ascoli's Theorem. Consequently, the sequence $\{\varphi_i^{\g_n}\}_{n\in\N}$ admits a subsequence that converges uniformly on an open neighborhood $U$ of $X$ which is relatively compact in $V_{0}$. By Hurwitz's Theorem, so will the associated derivatives subsequence. Moreover, in any analytic family, the function $\g\mapsto\varphi_j^\g(x)$, $\g\in\Ga$, is holomorphic for each $x\in V_{0}$ and $j\in\N$. Thus, the sequence $\{\varphi_i^{\g_n}\}_{n\in\N}$ converges pointwise to $\varphi_i^{\g_0}$ on $V_{0}$. Hence $\{\varphi_i^{\g_n}\}_{n\in\N}$ has a subsequence that converges uniformly on $U$ to $\varphi_i^{\g_0}$. The derivatives subsequence $\{(\varphi_i^{\g_n})'\}_{n\in\N}$ will then converge uniformly to $(\varphi_i^{\g_0})'$ on $U$, and in particular on $X$. This is a contradiction.  
\endpf

\

In light of Theorem~\ref{l1042903t}, it is relevant to ask whether 
an analytic family is continuous when $\CIFS(X)$ is endowed with 
the $\lambda$-topology. 

\begin{rem}\label{l1042903r}
An analytic family $\{\Phi^\g\}_{\g\in\Ga}$ regarded as the function 
\[ \begin{array}{ccccc} \Phi & : & \Ga & \rightarrow &\CIFS(X) \\
                             &   & \g & \mapsto     & \Phi^\g
   \end{array} 
\] 
is generally not continuous when $\CIFS(X)$ is endowed with the $\lambda$-topology.
\end{rem}
 
Indeed, if this were the case then the finiteness parameter function $\g\mapsto\theta_{\Phi^\g}$ would be locally constant according to Lemma~5.4 in~\cite{RU}. Since $\Ga$ is connected, that function would in fact be constant. However, Example~6 in section~8 of~\cite{RU} presents an analytic family whose finiteness parameter varies. 

\

An immediate consequence of Theorem~\ref{l1042903t} is the following.

\begin{cor}
Within an analytic family $\{\Phi^\g\}_{\g\in\Ga}$, 
the pressure function 
$(\g,t)\mapsto P(\g,t):=P(\Phi^\g,t)$,
is lower semi-continuous when $\CIFS(X)$ is endowed 
with the pointwise topology. Moreover, it is continuous 
on the set $\cup_{\g\in\Ga}
\{\g\}\times Fin(\Phi^\g)^c$.
\end{cor}

{\sl Proof.} 
It follows immediately from Theorem~\ref{l1042903t} that the map $(\g,t)\mapsto(\Phi^\g,t)$ is continuous, while Lemma~\ref{l1042903c} asserts that the map $(\Phi,t)\mapsto P(\Phi,t)$ is lower semi-continuous. The composition of these two functions, the pressure function $(\g,t)\mapsto P(\g,t):=P(\Phi^\g,t)$, is thus lower semi-continuous. Observe also that $P(\g,t)=\infty$ for any $t\notin Fin(\Phi^\g)$ and 
every $\g\in\Ga$. Thus, the pressure function is upper semi-continuous
at any point of the set $\cup_{\g\in\Ga}
\{\g\}\times Fin(\Phi^\g)^c$.
\endpf

\

\brem\label{rnc}
Note that $(\g,t)\mapsto P(\g,t)$
may be discontinuous at some points 
in $\cup_{\g\in\Ga}
\{\g\}\times Fin(\Phi^\g)$. For instance,
Examples~1, 2 and~3 in section~8 of~\cite{RU}
show that the function $\g\mapsto P(\g,t)$ 
may be discontinuous at $t=\theta_{\Phi^\g}$ if $\Phi^\g\notin\CFR(X)$. 
In all three examples, this function is discontinuous on the set
$\{(\g,t)\in\C\times\R:|\g|=1,t=1\}$.
\erem

\

We further have the following result when $I$ is finite.

\begin{thm}\label{subharm}
Let $I$ be a finite alphabet. If $\{\Phi^\g\}_{\g\in\Ga}$ 
is an analytic family in $\HIFS(X,I)$, then for every $t\ge 0$
the pressure function $\g\mapsto P(\g,t)$ 
is a continuous, plurisubharmonic function. 
\end{thm}

{\sl Proof.}  
Fix $t\geq 0$. Since $I$ is finite, each $\Phi^\g$ is regular. 
According to the proof of Theorem~6.3 in~\cite{RU}, 
fixing any $\g_{0}\in\Ga$, 
 we have that 
\[ P(\g,t) = \sup_{\mu\in M(\g_{0})} \left\{ h_{\mu}(\sigma) 
+ t\int_{I^{\infty}}\zeta_{\Phi^\g}(\omega) d\mu(\omega) \right\}, \]
and the function $\g\mapsto\int_{I^{\infty}} 
\zeta_{\Phi^\g}(\omega) d\mu(\omega)$ is finite, positive and pluriharmonic for every 
$\mu\in M(\g_{0})$. Thus, $P(\g,t)$ is the supremum of a family of 
pluriharmonic functions that are uniformly bounded from above by $\log|I|$. 
Therefore $P(\g,t)$ is a plurisubharmonic function uniformly 
bounded from above by $\log|I|$.
\endpf

\

In the infinite case, analytic families that are 
continuous with respect to the $\lambda$-topology on $\HIFS(X)$ enjoy the following property.

\bthm\label{psubharm}
Let $t\ge 0$ and $\{\Phi^\g\}_{\g\in\Ga}$ 
be an analytic family in $\HIFS(X)$ such that 
$\g\mapsto\Phi^{\g}\in\HIFS(X)$ is continuous with 
respect to the $\lambda$-topology. 
Suppose that there exists $\g_{0}\in\Ga$ 
such that $P(\g_{0},t)<\infty$. Then the function  
$\g\mapsto P(\g,t)$, $\g\in\Ga$, 
is continuous and plurisubharmonic. In particular, 
this function satisfies the Maximum Principle.    
\end{thm}

{\sl Proof.} 
Using Theorem~5.7 in~\cite{RU}, $\g\mapsto P(\g,t)$ is continuous. 
Moreover, using Theorem~2.1.5 in~\cite{gdms},  
the Variational Principle and Fubini's Theorem, it follows that the function 
$\g\mapsto P(\g,t)$ is a supremum of a family of 
pluriharmonic functions. Hence it is plurisubharmonic.   
\endpf

\

\subsection{Real Analyticity of Pressure}

Since all the forthcoming results rely upon the theorem below and since 
its proof is somewhat lengthy we single it out in this subsection. 
Moreover, in order to allege notation 
we will write $\pi_\g$ for $\pi_{\Phi^\g}$ and
$P(\g,t)$ for $P(\Phi^\g,t)$.

\bthm\label{prealanalconj}
Let $\{\Phi^\g\}_{\g\in\Ga}$
be an analytic family in $\HIFS(X)$ with the following property.
\begin{itemize}
\item 
For every $\g_0\in\Ga$ there exist a neighborhood
$\Ga_0$ of $\g_0$ and $B\ge 0$ such that 
$$
\lt|\frac{(\varphi_{\omega_1}^\g)'(\pi_\g(\sigma(\omega)))}
{(\varphi_{\omega_1}^{\g_0})'(\pi_{\g_0}(\sigma(\omega)))}\rt|\le B 
$$ 
for every $\omega\in I^\infty$ and all $\g\in\Ga_0$.
\end{itemize}
If $t\ge 0$ and $P(\g_1,t)<\infty$ for some $\g_1\in\Ga$, then
$P(\g,t)<\infty$ for all $\g\in \Ga$ and the function $\g\mapsto 
P(\g,t)$, $\g\in\Ga$, is real-analytic.
\ethm

\

\ni In the case $\dim(\Ga)=1$ the proof of this theorem is contained in
the proof of Theorem~4.2 in \cite{U}, where one assumption is weaker and one
is somewhat stronger than in the present paper. Here, we actually provide a
different, simplified proof with $\dim(\Ga)$ an arbitrary finite positive integer.
Even this simplified proof requires substantial preparation. We begin 
as follows. Denote by $\N_0$ the set $\N\cup\{0\}$. Let $d\ge 1$. For every $\g_0=(\g_1^0,\g_2^0,\ldots,\g_d^0)\in\C^d$, denote by $\C^dF(\g_0)$ the 
real vector space of all formal power series 
$$
f(\g_1,\g_2,\ld,\g_d)
=\sum_{\a\in\N_0^d}f_\a(\g_1-\g_1^0)^{\a_1}(\g_2-\g_2^0)^{\a_2}\ld
                       (\g_d-\g_d^0)^{\a_d}, \ (\g_1,\g_2,\ld,\g_d)\in\C^d,
$$
with complex coefficients $f_\a$.
Similarly, for every $\g_0\in\R^d$, denote by $\R^dF(\g_0)$ 
the real vector space of all formal power series 
$$
g(\g_1,\g_2,\ld,\g_d)
=\sum_{\a\in\N_0^d}g_\a(\g_1-\g_1^0)^{\a_1}(\g_2-\g_2^0)^{\a_2}\ld
                       (\g_d-\g_d^0)^{\a_d}, \ (\g_1,\g_2,\ld,\g_d)\in\R^d,
$$
with real coefficients $g_\a$. Moreover, for every $\g_0\in\R^d$, let 
$U_d^{\g_0}:\R^dF(\g_0)\to\C^dF(\g_0)$ be the linear operator which ascribes 
to each formal power series 
$\sum_{\a\in\N_0^d}g_\a(\g_1-\g_1^0)^{\a_1}(\g_2-\g_2^0)^{\a_2}\ld (\g_d-\g_d^0)^{\a_d}$
in $\R^dF(\g_0)$ the formal power series 
$\sum_{\a\in\N_0^d}g_\a(\g_1-\g_1^0)^{\a_1}(\g_2-\g_2^0)^{\a_2}\ld (\g_d-\g_d^0)^{\a_d}$
in $\C^dF(\g_0)$. Furthermore, for every $\g_0\in\C^d$, set
$$
I(\g_0)=\bigl(\re(\g_1^0),\im(\g_1^0),\re(\g_2^0),\im(\g_2^0),\ld,\re(\g_d^0),\im(\g_d^0)\bigr)\in\R^{2d},
$$ 
and let $\re_d^{\g_0}:\C^dF(\g_0)\to\R^{2d}F(I(\g_0))$ be the 
operator which ascribes to every formal power series 
$$
f(\g_1,\g_2,\ld,\g_d)
=\sum_{\a\in\N_0^d}f_\a(\g_1-\g_1^0)^{\a_1}(\g_2-\g_2^0)^{\a_2}\ld
                       (\g_d-\g_d^0)^{\a_d}, \ (\g_1,\g_2,\ld,\g_d)\in\C^d,
$$
in $\C^dF(\g_0)$ the formal power series $\re_d^{\g_0}(f)\in\R^{2d}F(I(\g_0))$ 
which is given by
$$
\re_d^{\g_0}(f)(x_1,y_1,x_2,y_2,\ld,x_d,y_d)
=\sum_{\b\in\N_0^{2d}}c_\b\prod_{j=1}^d\(x_j-\re(\g_j^0)\)^{\b_j^{(1)}}\(y_j-\im(\g_j^0)\)^{\b_j^{(2)}},
$$
where each $\b\in\N_0^{2d}$ is written in the form
$\(\b_1^{(1)},\b_1^{(2)},\b_2^{(1)},\b_2^{(2)},\ld,\b_d^{(1)},\b_d^{(2)}\)$ and 
$$
c_\b=\re\lt(f_{\hat\b}\prod_{j=1}^d
      {\b_j^{(1)}+\b_j^{(2)}\choose\b_j^{(1)}}i^{\b_j^{(2)}}\rt)
$$
with $\hat\b=
\(\b_1^{(1)}+\b_1^{(2)},\b_2^{(1)}+\b_2^{(2)},\ld,\b_d^{(1)}+\b_d^{(2)}\)\in\N_0^d$.
Clearly, $f\mapsto\re_d^{\g_0}(f)$ is a linear operator from the real vector space
$\C^dF(\g_0)$ to the real vector space $\R^{2d}F(I(\g_0))$. Now, fix $R>0$
and consider $\Hol_b^d(D_d(\g_0,R))$, the real vector space of all formal power
series in $\C^dF(\g_0)$ that converge uniformly on the $d$-dimensional disk $D_d(\g_0,R)$. These series have
finite supremum norms over $D_d(\g_0,R)$. Moreover, $\Hol_b^d(D_d(\g_0,R))$ endowed with the 
supremum norm over $D_d(\g_0,R)$ is a normed vector space. Let
$$
\Ra_d^{\g_0}=U_{2d}^{I(\g_0)}\circ\re_d^{\g_0}:\C^dF(\g_0)\to \C^{\,2d}F(I(\g_0)).
$$
We will need the following fact. This fact was proved in \cite{myu} as Lemma~8.1.

\

\blem\label{l1myu9.1}
For every $\g_0\in\C^d$ and every $R>0$ we have 
$$
\Ra_d^{\g_0}(\Hol_b^d(D_d(\g_0,R)))\sbt \Hol_b^{2d}(D_{2d}(I(\g_0),R/4)).
$$
Moreover, the linear operator $\Ra_d^{\g_0}:
\Hol_b^d(D_{d}(\g_0,R))\to\Hol_b^{2d}(D_{2d}(I(\g_0),R/4))$ is bounded and its norm 
is bounded above by $4^d$. Furthermore, for every $f\in\Hol_b^d(D_{d}(\g_0,R))$, 
$$
\Ra_d^{\g_0}(f)|_{D_{d}(\g_0,R/4)}=\re(f)|_{D_{d}(\g_0,R/4)},
$$
where $\re(f)$ is the real part of the complex-valued function $f$ and 
$D_{d}(\g_0,R/4)$ is identified with $D_{2d}(I(\g_0),R/4)\cap\R^{2d}$ as 
$\C^d$ is embedded in $\C^{\,2d}$ by the formula $(x_1+iy_1,x_2+iy_2,\ld,x_d+iy_d)
\mapsto (x_1,y_1,x_2,y_2,\ld,x_d,y_d)$.
\elem

\

\ni Passing to other ingredients needed in our proof of Theorem~\ref{prealanalconj},
we recall the space of H\"older continuous functions. Let $\a>0$. We say that a 
function $g:I^\infty\to\C$ is 
$\a$-H\"older continuous if 
$$
v_\a(g):=\sup_{n\in\N}\{v_{\a,n}(g)\}<\infty,
$$
where  
$$
v_{\a,n}(g)
=\sup\{|g(\om)-g(\tau)|e^{\a n}:\om,\tau\in I^\infty
\text{ and } |\om\wedge \tau|\ge n\}.
$$
Let $\K_\a$ be the set of all complex-valued $\a$-H\"older continuous (but not 
necessarily bounded) functions defined on $I^\infty$. 
Set 
$$
\K_{\a }^{S}
:=\lt\{u\in \K_\a:\sum_{i\in I}\exp\(\sup_{\om\in[i]}\re(u(\om))\)<\infty\rt\}.
$$
Members of $\K_\a^S$ are called $\a$-H\"older summable potentials. 
If a function $g:I^\infty\to\C$ has finite supremum norm $\|g\|_\infty$ and 
finite $v_\a(g)$, then $g$ is said to belong to the space $\Ha_\a$ 
of all bounded $\a$-H\"older continuous functions
defined on $I^\infty$. $\Ha_\a$ becomes a Banach space when equipped with the norm 
$\|\cdot\|_\a$ given by the formula
$$
\|g\|_\a=\|g\|_\infty + v_\a(g).
$$

\

\blem\label{m0531080}
For every $g, h\in\Ha_\a$ we have 
$\|gh\|_\a \le 3\|g\|_\a\|h\|_\a.$
\elem

{\sl Proof.} Indeed, $\|gh\|_\infty\le\|g\|_\infty\|h\|_\infty\le\|g\|_\a\|h\|_\a$. 
Moreover, for every $\om, \tau\in I^\infty$ such that $|\om\wedge\tau|\ge 1$ we have
\begin{eqnarray*} 
|(gh)(\om)-(gh)(\tau)|
&\le&|g(\om)|\,|h(\om)-h(\tau)|+|h(\tau)|\,|g(\om)-g(\tau)| \\
&\le&\|g\|_\infty \, v_\a(h)e^{-\a|\om\wedge\tau|}+\|h\|_\infty \, v_\a(g)e^{-\a|\om\wedge\tau|} \\
&\le&2\|g\|_\a\|h\|_\a e^{-\a|\om\wedge\tau|}, 
\end{eqnarray*}
that is, $v_\a(gh)\le 2\|g\|_\a\|h\|_\a$. Consequently, $\|gh\|_\a\le 3\|g\|_\a\|h\|_\a$.
\endpf

\

In particular, this lemma shows that $\Ha_\a$ becomes a Banach algebra 
if one replaces $\|\cdot\|_\a$ by an equivalent norm $\|\cdot\|_\a'$ 
such that $\|gh\|_\a'\le\|g\|_\a'\|h\|_\a'$. The factors ``3'' and ``9''
in the following lemma would then ``disappear''.

\

\blem\label{m0531081}
For every $g\in\Ha_\a$, 
$$
\|e^g\|_\a\le e^{3\|g\|_\a} 
\hspace{1cm} \mbox{ and } \hspace{1cm} 
\|e^g-1\|_\a\le 3\|g\|_\a e^{3\|g\|_\a},
$$
and for every $g,h\in\Ha_\a$, 
$$
\|e^h-e^g\|_\a\le 9e^{3(\|g\|_\a+\|h-g\|_\a)}\|h-g\|_\a.
$$
In particular, the map $g\in\Ha_\a\mapsto e^g\in\Ha_\a$ is continuous. 
\elem

{\sl Proof.} 
Since $e^g=\sum_{n=0}^\infty g^n/n!$, it follows from Lemma~\ref{m0531080}
that $\|e^g\|_\a\le e^{3\|g\|_\a}$. Moreover,   
\begin{eqnarray*}
\|e^g-1\|_\a&=&\bigl\|\sum_{n=1}^\infty \frac{g^n}{n!}\bigr\|_\a
            =\bigl\|g\sum_{n=0}^\infty \frac{g^n}{(n+1)!}\bigr\|_\a \\
            &\le& 3\|g\|_\a \bigl\|\sum_{n=0}^\infty \frac{g^n}{(n+1)!}\bigr\|_\a
            \le 3\|g\|_\a \sum_{n=0}^\infty \frac{\|g^n\|_\a}{(n+1)!} \\
            &\le& 3\|g\|_\a \sum_{n=0}^\infty \frac{(3\|g\|_\a)^n}{(n+1)!}
            \le 3\|g\|_\a \sum_{n=0}^\infty \frac{(3\|g\|_\a)^n}{n!}
            = 3\|g\|_\a e^{3\|g\|_\a}.
\end{eqnarray*}
Hence, using Lemma~\ref{m0531080} and the above two estimates, we deduce that for every $g,h\in\Ha_\a$, 
$$
\|e^h-e^g\|_\a=\|e^g(e^{h-g}-1)\|_\a
\le 3\|e^g\|_\a\|e^{h-g}-1\|_\a
\le 9e^{3(\|g\|_\a+\|h-g\|_\a)}\|h-g\|_\a.
$$
\endpf

\

Now, for every $u\in \K_\a^S$ and $g\in\Ha_\a$ define 
$$
\L_ug(\om)=\sum_{i\in I}e^{u(i\om)}g(i\om).
$$
Easy calculations 
show that $\L_ug\in\Ha_\a$ and that the operator $\L_u:\Ha_\a\to\Ha_\a$  
is linear and bounded. We hence have the following.

\

\blem\label{l1rs253}
If $u\in \K_\a^S$, then 
$\L_u\in L(\Ha_\a)$, where $L(\Ha_\a)$ denotes the Banach space of all
bounded, linear operators on $\Ha_\a$, endowed with the operator norm. 
\elem

\ni We immediately obtain the following.

\

\bcor\label{l1rs2531}
Let $\Phi\in\HIFS(X)$ and $s\in(0,1)$ be a contraction ratio for 
$\Phi$. Let $0<\a\le -\log s$. If $P(\Phi,t)<\infty$ and $g\in\Ha_\a$, then the potential $u(\om)=g(\om)+
t\log|\varphi_{\om_1}'(\pi_\Phi(\sg(\om)))|$ belongs to $\K_\a^S$ and the 
corresponding linear operator $\L_{u}:\Ha_\a\to\Ha_\a$ is bounded.
\ecor

{\sl Proof.} We simply need to prove that $u\in\K_\a^S$.  
To do this, observe that by Remark~\ref{rsbdp} and Lemma~\ref{lsbdp}, 
for all $\om,\tau\in I^\infty$ with $|\om\wedge\tau|\ge 1$
we have 
$$
\aligned
\left|\log|\varphi_{\om_1}'(\pi_\Phi(\sg(\om)))|-
      \log|\varphi_{\tau_1}'(\pi_\Phi(\sg(\tau)))|\right|
&\leq L(1-s)^{-1}|\pi_\Phi(\sg(\om))-\pi_\Phi(\sg(\tau))|\\ 
&\leq L(1-s)^{-1} e^{-\a(|\om\wedge\tau|-1)}\mbox{diam}(X). 
\endaligned
$$
Thus, the function  
$\om\mapsto\log|\varphi_{\om_1}'(\pi_\Phi(\sg(\om)))|$ is in $\K_\a$, as $v_\a(\om\mapsto\log|\varphi_{\om_1}'(\pi_\Phi(\sg(\om)))|)
\le L(1-s)^{-1}e^\a\mbox{diam}(X)$. Since $g\in\K_\a$ as well, 
we conclude that $u\in\K_\a$. 
Moreover, 
\begin{eqnarray*}
\sum_{i\in I}\exp\(\sup_{\om\in[i]}\re(u(\om))\)
&=&\sum_{i\in I}\exp\(\sup_{\om\in[i]}\re\(g(\om)+t\log|\varphi_i'(\pi_\Phi(\sg(\om)))|\)\) \\
&\le&\sum_{i\in I}\exp\(\sup_{\om\in[i]}\re(g(\om))\)\|\varphi_i'\|^t \\
&\le&e^{\|g\|_\infty}\sum_{i\in I}\|\varphi_i'\|^t<\infty 
\end{eqnarray*}
since $P(\Phi,t)<\infty$. Thus, $u$ is summable.
\endpf

\

\ni We shall now prove the following result. Note that the finite-dimensional 
complex manifold $\^\Ga$ in this result is completely general. It need not 
be the same and is generally not the same as the manifold $\Ga$ from 
Theorem~\ref{prealanalconj}. This important 
and necessary distinction will become clearer in Lemma~\ref{l1rs259} , where 
$\Ga$ will be the manifold from Theorem~\ref{prealanalconj} and is a subset 
of $\C^d$, whereas $\^\Ga$ will be $\C^{\,2d}$. Thus, in our case,
$\Ga$ can be naturally embedded in $\^\Ga$ via the embedding of $\C^d$ in 
$\C^{\,2d}$ given by the formula $(x_1+iy_1,x_2+iy_2,\ld,x_d+iy_d)
\mapsto (x_1,y_1,x_2,y_2,\ld,x_d,y_d)$.

\

\blem\label{l2rs253}
Let $G$ be an open subset of 
a finite-dimensional complex manifold $\^\Ga$.
Let $\^\g_0\in G$ and $\Phi^{\^\g_0}\in\HIFS(X)$.
Let also $t\ge 0$ be such that $P(\^\g_0,t)<\infty$, and 
$s\in(0,1)$ a contraction ratio for $\Phi^{\^\g_0}$. 
Set $0<\a\le -\log s$. Let 
${\^\g}\mapsto g_{\^\g}\in\Ha_\a$ be a function with the following properties:
\begin{itemize}
\item[(a)] $\sup\{\|g_{\^\g}\|_\a:\^\g\in G\}<\infty$; and
\item[(b)] for every $\om\in I^\infty$ the function ${\^\g}\mapsto g_{\^\g}(\om)$ 
from $G$ to $\C$ is holomorphic.
\end{itemize}
Set $\L_{\^\g}:=\L_{u_{\^\g}}$, where 
$u_{\^\g}(\om)=g_{\^\g}(\om)+t\log|(\varphi_{\om_1}^{\^\g_0})'(\pi_{\^\g_0}(\sg(\om)))|$. 
Then the map ${\^\g}\mapsto \L_{\^\g}\in L(\Ha_\a)$ is holomorphic.
\elem

\

In order to prove this result, we need a few preparatory lemmas. 
Moreover, in order to allege notation, we use $\g$ instead of ${\^\g}$ 
from this point on until the end of the proof of Lemma~\ref{l2rs253}. 
Also, observe that 
\[ \L_\g f(\om)=\sum_{i\in I}e^{u_\g(i\om)} f(i\om)
               =\sum_{i\in I}e^{g_\g(i\om)}|(\varphi_i^{\g_0})'(\pi_{\g_0}(\om))|^t f(i\om) \]
for every $\om\in I^\infty$. For every $i\in I$ and $\g\in G$, 
we define the linear operator $\L_{\g,i}:\Ha_\a\to\Ha_\a$ as 
$$
\L_{\g,i}f(\om)=e^{u_\g(i\om)} f(i\om)
               =e^{g_\g(i\om)}|(\varphi_i^{\g_0})'(\pi_{\g_0}(\om))|^t f(i\om).
$$
Finally, we let $H:=\sup\{\|g_{\overline{\g}}\|_\a:{\overline{\g}}\in G\}$.

\

%
%
%
%

\blem\label{m0531083}
Under the assumptions of Lemma~\ref{l2rs253},
the function $\g\mapsto g_\g$ from $G$ to $\Ha_\a$ is continuous.
\elem

{\sl Proof.} Applying Cauchy's formula to the holomorphic functions $\g\mapsto g_\g(\om)$, $\om\in I^\infty$,
we obtain that for each $\zeta\in G$ and each $\e>0$, there exists a $\delta_1>0$ such that 
for each $\g\in G$ with $|\g-\zeta|<\delta_1$ and for all $\om\in I^\infty$, 
$$
|g_\g(\om)-g_\zeta(\om)|<\e H.
$$
Hence for each $\g\in G$ with $|\g-\zeta|<\delta_1$, 
\begin{equation}
\label{sumiremeq1}
\|g_\g-g_\zeta\|_\infty<\e H.
\end{equation}

Moreover, for every $\g,\xi,\zeta\in G$ and for all $\om,\tau\in I^\infty$
such that $|\om\wedge\tau|\ge 1$,
\begin{eqnarray*}
\bigl|\(g_\g(\om)-g_\zeta(\om)\)-\(g_\g(\tau)-g_\zeta(\tau)\)\bigr|
&=&\bigl|\(g_\g(\om)-g_\g(\tau)\)-\(g_\zeta(\om)-g_\zeta(\tau)\)\bigr| 
\end{eqnarray*}
and 
$$
|g_\xi(\om)-g_\xi(\tau)|
                      \le H e^{-\a|\om\wedge\tau|}.
$$
Applying Cauchy's formula to the holomorphic functions 
$\g\mapsto g_\g(\om)-g_\gamma(\tau)$, $\om,\tau\in I^\infty$,
and using the previous two relations and condition (a), 
we obtain that for each $\zeta\in G$ and each $\e>0$, there 
exists a $\delta_2>0$ such that 
for each $\g\in G$ with $|\g-\zeta|<\delta_2$ and 
for all $\om, \tau\in I^\infty$ such that $|\om\wedge\tau|\ge 1$, 
$$
\bigl|\(g_\g(\om)-g_\zeta(\om)\)-\(g_\g(\tau)-g_\zeta(\tau)\)\bigr|
\leq \e H e^{-\a|\om\wedge\tau|}. 
$$ 
Hence for each $\g\in G$ with $|\g-\zeta|<\delta_2$, 
\begin{equation}
\label{sumiremeq2}
v_\a(g_\g-g_\zeta)\leq \e H.
\end{equation}
From (\ref{sumiremeq1}) and (\ref{sumiremeq2}), 
it follows that for each $\g\in G$ with $|\g-\zeta|<\min\{\delta_1,\delta_2\}$, 
$$
\|g_\g-g_\zeta\|_\a\le 2\e H.
$$ 
Therefore, the function $\g\mapsto g_\g$ is continuous. 
\endpf

\

\blem\label{m0531084} 
$\||(\varphi_i^{\g_0})'\circ\pi_{\g_0}|^t\|_\a
\le C \|(\varphi_i^{\g_0})'\|^t$, where 
$C=1+KLt\diam(X)$.
\elem

{\sl Proof.} In view of conditions (iv) and (iv') defining CIFSs 
in section~\ref{secprelifs}, we get for all $\om,\tau\in I^\infty$ 
with $|\om\wedge\tau|\ge 1$ that 
\begin{eqnarray*}
\big||(\varphi_i^{\g_0})'(\pi_{\g_0}(\om))|^t-|(\varphi_i^{\g_0})'(\pi_{\g_0}(\tau))|^t\big|
&\le& t\|(\varphi_i^{\g_0})'\|^{t-1}
\big||(\varphi_i^{\g_0})'(\pi_{\g_0}(\om))|-|(\varphi_i^{\g_0})'(\pi_{\g_0}(\tau))|\big| \\
&\le& t\|(\varphi_i^{\g_0})'\|^{t-1}
L\|[(\varphi_i^{\g_0})']^{-1}\|_V^{-1}|\pi_{\g_0}(\om)-\pi_{\g_0}(\tau)| \\
&\le& t\|(\varphi_i^{\g_0})'\|^{t-1}
LK\|(\varphi_i^{\g_0})'\|s^{|\om\wedge\tau|}\mbox{diam}(X) \\
&\le& KLt\|(\varphi_i^{\g_0})'\|^t e^{-\a|\om\wedge\tau|}\mbox{diam}(X) \\
&\le& KLt\mbox{diam}(X)\|(\varphi_i^{\g_0})'\|^t e^{-\a|\om\wedge\tau|}.
\end{eqnarray*}
Therefore $v_\a(|(\varphi_i^{\g_0})'\circ\pi_{\g_0}|^t)
\le KLt\mbox{diam}(X)\|(\varphi_i^{\g_0})'\|^t$. 
It is then straightforward to see that $\||(\varphi_i^{\g_0})'\circ\pi_{\g_0}|^t\|_\a
\le (1+KLt\diam(X)) \|(\varphi_i^{\g_0})'\|^t$.
\endpf

\

\blem\label{m0531085}
Under the assumptions of Lemma~\ref{l2rs253},
for each $i\in I$ the function $\g\mapsto\L_{\g,i}$ 
from $G$ to $L(\Ha_\a)$ is continuous.
Moreover, $\|\L_{\g,i}\|\le 9e^{3H}C \|(\varphi_i^{\g_0})'\|^t$.
\elem

{\sl Proof.} From the expression of $\L_{\g,i}f$ and 
Lemmas~\ref{m0531080},~\ref{m0531081}  
and~\ref{m0531084}, we deduce that for every $\g,\zeta\in G$ 
and for every $f\in\Ha_\a$,
\begin{eqnarray*}
\|\L_{\g,i}f-\L_{\zeta,i}f\|_\a
&\le& 9\|e^{g_\g}-e^{g_\zeta}\|_\a \||(\varphi_i^{\g_0})'\circ\pi_{\g_0}|^t\|_\a \|f\|_\a \\
&\le& 9\cdot 9e^{3(\|g_\g\|_\a+\|g_\g-g_\zeta\|_\a)}\|g_\g-g_\zeta\|_\a \cdot C \|(\varphi_i^{\g_0})'\|^t \cdot \|f\|_\a \\
&\le& 81C \|(\varphi_i^{\g_0})'\|^t e^{3(\|g_\g\|_\a+\|g_\g-g_\zeta\|_\a)}\|g_\g-g_\zeta\|_\a \|f\|_\a.
\end{eqnarray*}
Therefore $\|\L_{\g,i}-\L_{\zeta,i}\|\le 81C \|(\varphi_i^{\g_0})'\|^t e^{3(\|g_\g\|_\a+\|g_\g-g_\zeta\|_\a)}\|g_\g-g_\zeta\|_\a$
and the continuity of $\g\mapsto\L_{\g,i}$ follows from that 
of $\g\mapsto g_\g$, that is, from Lemma~\ref{m0531083}.

Furthermore, from the expression of $\L_{\g,i}f$ and 
Lemmas~\ref{m0531080},~\ref{m0531081}  and~\ref{m0531084}, 
we obtain that 
\begin{eqnarray*}
\|\L_{\g,i}f\|_\a
&\le& 9\|e^{g_\g}\|_\a \||(\varphi_i^{\g_0})'\circ\pi_{\g_0}|^t\|_\a \|f\|_\a \\
&\le& 9e^{3\|g_\g\|_\a}\cdot C \|(\varphi_i^{\g_0})'\|^t \cdot \|f\|_\a \\
&\le& 9e^{3H}C \|(\varphi_i^{\g_0})'\|^t \|f\|_\a.
\end{eqnarray*}
Thus, $\|\L_{\g,i}\|\le 9e^{3H}C \|(\varphi_i^{\g_0})'\|^t$.
\endpf

\

%

We now prove the holomorphicity of that function.

\ 

\blem\label{m0531087}
Under the assumptions of Lemma~\ref{l2rs253}, for each $i\in I$
the function $\g\mapsto\L_{\g,i}$ from $G$ to $L(\Ha_\a)$ is holomorphic. 
\elem

{\sl Proof.} Indeed, from Lemma~\ref{m0531085} and Osgood's 
Lemma for (Banach space)-valued functions (see the 
argument on page 3 of~\cite{G}), 
we may assume that $\dim_{\C}\tilde{\Ga}=1$. 
Fix $i\in I$. Let $f\in\Ha_\a$, $\om\in I^\infty$ and 
let $T$ be any triangle in $G$. Since the function $\g\mapsto\L_{\g,i}f(\om)$ 
from $G$ to $\C$ is holomorphic by condition (b) of 
Lemma~\ref{l2rs253}, we have that 
\[ \int_{\partial T}\L_{\g,i}f(\om) \, d\g=0. \]
Therefore  
\[ \left(\int_{\partial T}\L_{\g,i} \, d\g\right)(f)(\om)=\int_{\partial T}\L_{\g,i}f(\om) \, d\g=0. \]
Since this holds for all $\om \in I^\infty$ and for all $f\in\Ha_\a$, we deduce that 
\[ \int_{\partial T}\L_{\g,i} \, d\g=0. \]
As the triangle $T$ was chosen arbitrarily in $G$, we conclude from Morera's Theorem that 
the function $\g\mapsto\L_{\g,i}$ 
from $G$ to $L(\Ha_\a)$ is holomorphic. 
\endpf 

\

{\sl Proof of Lemma~\ref{l2rs253}.} Since $\L_\g=\sum_{i\in I}\L_{\g,i}$, 
since all the functions $\g\mapsto \L_{\g,i}$ are holomorphic
according to Lemma~\ref{m0531087}, since $\|\L_{\g,i}\|\le 9e^{3H}C \|(\varphi_i^{\g_0})'\|^t$
for every $i\in I$ according to Lemma~\ref{m0531085}, since the series $\sum_{i\in I}\|(\varphi_i^{\g_0})'\|^t$
converges as $P(\g_0,t)<\infty$, and since the sum of a uniformly convergent sequence of holomorphic functions is holomorphic, we deduce that the function $\g\mapsto \L_\g$ 
from $G$ to $L(\Ha_\a)$ is holomorphic. 
\endpf
 
\

We continue our preparation by looking at the coding map.

\

\blem 
\label{l:pigomegaanal}
Let $\{\Phi^\g\}_{\g\in\Ga}$ be an analytic 
family in $\HIFS(X)$. For any fixed $\om\in I^\infty$, 
the function $\g\mapsto\pi_\g(\om)$
is analytic and the family $\{\g\mapsto\pi_\g(\om)\}_{\om\in I^\infty}$
is equicontinuous on $\Ga$. 
\elem
{\sl Proof.} 
Let $x\in X$. Fix for a moment $\om\in I^\infty$. 
Then $\pi_\g(\om)=
\lim_{n\rightarrow\infty}\varphi_{\om_1}^\g\circ 
\cdots\circ\varphi_{\omega_n}^\g(x)$. 
Hence, $\g\mapsto\pi_\g(\om)$ is a 
limit of a sequence of uniformly bounded holomorphic maps defined on 
$\Ga$. Since a family of uniformly bounded holomorphic maps is 
normal, it follows that $\g\mapsto\pi_\g(\omega)$ is 
holomorphic. Furthermore, 
$\{\g\mapsto\pi_\g(\omega)\}_{\omega\in I^\infty}$ 
is uniformly bounded. Hence, it is normal and thereby equicontinuous. 
\endpf 

\

We now introduce a function $\kappa$ and describe some of its properties.

\

\blem\label{l2rs2532}
In the setting of Theorem~\ref{prealanalconj}, 
fix $\g_0\in\Ga$ and $R_0>0$ so that
$D_d(\g_0,R_0)\sbt \Ga_0$. For every $\om\in I^\infty$ define
$$
\ka_\om(\g)={(\varphi_{\om_1}^\g)'(\pi_\g(\sg(\om)))\over 
             (\varphi_{\om_1}^{\g_0})'(\pi_{\g_0}(\sg(\om)))},
\,\,\,\g\in D_d(\g_0,R_0).
$$
Then we have that
$$
|\ka_\om(\g)-1|\le {1\over 4}
$$
for all $\om\in I^\infty$ and all $\g\in D_d(\g_0,R_0/8B)$,
and thus for all $\om\in I^\infty$ there exists a unique holomorphic 
branch $\log\ka_\om$ of the logarithm of $\ka_\om$ defined on $D_d(\g_0,R_0/8B)$
and sending $\g_0$ to $0$. Furthermore, there exists a constant $M>0$ such that
$$
|\log\ka_\om(\g)|\le M
$$
for all $\om\in I^\infty$ and all $\g\in D_d(\g_0,R_0/8B)$. 
\elem

{\sl Proof.} According to Condition (1) from Theorem~\ref{prealanalconj} 
we have that
\begin{equation}\label{1h51}
|\ka_\om(\g)|\le B, 
\end{equation}
for all $\om\in I^\infty$ and all $\g\in D_d(\g_0,R_0)$.
Choosing an appropriate rotation to
change the system of coordinates we may assume without loss of generality that the
dimension $d=1$. Since for all $\g\in B(\g_0,R_0/2)$ and all $0<r\le R_0/2$,
$$
\ka_\om'(\g)={1\over 2\pi i}\int_{\bd B(\g,r)}{\ka_\om(\xi)\over (\xi-\g)^2}d\xi,
$$
we obtain from (\ref{1h51}) that
$$
|\ka_\om'(\g)|
\le {1\over 2\pi}\int_{\bd B(\g,r)}{B\over r^2}|d\xi|
={B\over r}.
$$
Since $\ka_\om(\g_0)=1$, we therefore get for all $\g\in B(\g_0,r)$ that
$$
|\ka_\om(\g)-1|
=|\ka_\om(\g)-\ka_\om(\g_0)|
=\lt|\int_{\g_0}^\g\ka_\om'(\xi)d\xi\rt|
\le \int_{\g_0}^\g|\ka_\om'(\xi)||d\xi|
\le {B\over r}|\g-\g_0|.
$$
Taking $r=R_0/2$, we deduce that $|\ka_\om(\g)-1|\le 1/4$ for all $\g\in D_d(\g_0,R_0/8B)$. \endpf

\

Our last auxiliary result is the following.

\

\blem\label{l1rs259}
In the setting of Theorem~\ref{prealanalconj}, 
fix $\g_0\in\Ga$ and $R_0>0$ so that
$D_d(\g_0,R_0)\sbt \Ga_0$. Let $s\in(0,1)$ 
be a contraction ratio for $\Phi^{\g_0}$. Set $0<\a<-\log s$
and choose $0<R<R_0/8B$ so small that $e^{-\a}$ is a contraction 
ratio for all $\g\in D_d(\g_0,R)$. Then the map 
${\^\g}\mapsto \L_{\^\g}\in L(\Ha_\a)$, ${\^\g}\in D_{2d}(I(\g_0),R/4)$,
 is holomorphic, where 
$\L_{\^\g}:=\L_{u_{\^\g}}$ with 
\[ u_{\^\g}(\om)=t\Ra_d^{\g_0}(\log\ka_\om)(\^\g)+t
               \log|(\varphi_{\om_1}^{\g_0})'(\pi_{\g_0}(\sg(\om)))|. \] 
\elem

{\sl Proof.} Let $\^\Ga=\C^{\,2d}$, $G=D_{2d}(I(\g_0),R/4)$, $\^\g_0=I(\g_0)$, $\Phi^{\^\g_0}=\Phi^{\g_0}$, and let $g_{\^\g}(\om)=t\Ra_d^{\g_0}(\log\kappa_\om)(\^\g)$. According to Lemma~\ref{l2rs253}, it suffices to prove that for all ${\^\g}\in D_{2d}(I(\g_0),R/4)$ the function $\om\mapsto\Ra_d^{\g_0}(\log\ka_\om)(\^\g)$ 
belongs to $\Ha_\a$, and that the function $\^\g\mapsto\Ra_d^{\g_0}(\log\ka)(\^\g)$ from $D_{2d}(I(\g_0),R/4)$ 
to $\Ha_\a$ satisfies conditions (a) and (b). Condition (b)  
is clearly satisfied. Moreover, it follows from 
Lemmas~\ref{l1myu9.1} and~\ref{l2rs2532} that 
\begin{equation}\label{1rs259}
\|\Ra_d^{\g_0}(\log\ka)(\^\g)\|_\infty\le 4^dM,
\end{equation}
where $M$ is the constant arising from Lemma~\ref{l2rs2532}.


\

Claim: There exists a constant $C_0>0$ such that 
for all $\om, \tau\in I^\infty$ and for all $\g\in D_d(\g_0,R)$, 
$$|\log\ka_\om(\g)-\log\ka_\tau(\g)|\leq C_0 e^{-\a|\om\wedge\tau|}.$$

In order to prove the above claim, observe that if $\om_1\neq\tau_1$, then 
$$|\log\ka_\om(\g)-\log\ka_\tau(\g)|
\leq |\log\ka_\om(\g)|+|\log\ka_\tau(\g)|
\leq 2M = 2M e^{-\a|\om\wedge\tau|}.$$
Thus, suppose $\om _{1}=\tau _{1}=i$ for some $i\in I$. Let $\beta_{\g,\om}:=\pi_\g(\sg(\om))$. Then  
\begin{equation}\label{mar2}
\aligned
|\ka_\om(\g)&-\ka_\tau(\g)|\\ 
& =\left|\frac{(\varphi_i^\g)'(\pi_\g(\sg(\om)))(\varphi_i^{\g_0})'(\pi_{\g_0}(\sg(\tau)))
-(\varphi_i^\g)'(\pi_\g(\sg(\tau)))(\varphi_i^{\g_0})'(\pi_{\g_0}(\sg(\om)))}
{(\varphi_i^{\g_0})'(\pi_{\g_0}(\sg(\om)))(\varphi_i^{\g_0})'(\pi_{\g_0}(\sg(\tau)))}\right|\\ 
& \leq K\frac{|(\varphi_i^{\g_0})'(\beta_{\g_0,\tau})| 
|(\varphi_i^\g)'(\beta_{\g,\om})-(\varphi_i^\g)'(\beta_{\g,\tau})|}
{|(\varphi_i^{\g_0})'(\beta_{\g_0,\tau})|^{2}}\\ 
&\ \ \hspace{1cm}+K\frac{|(\varphi_i^\g)'(\beta_{\g,\tau})| 
|(\varphi_i^{\g_0})'(\beta_{\g_0,\tau})-(\varphi_i^{\g_0})'(\beta_{\g_0,\om})|}
{|(\varphi_i^{\g_0})'(\beta_{\g_0,\tau})|^{2}}.
\endaligned 
\end{equation}
Let $r:=\mbox{dist}(\partial V_0,X)$, where $V_0$ denotes the neighborhood of $X$ from Definition 3.2. 
Let $\{B(z_j,r/2)\}_{j=1}^J$ be a finite covering of 
$X$ where $z_j\in X$ for each $j$. 
For each $1\leq j\leq J$, $i\in I$ and $\g\in D_d(\g_0,R)$, we define 
$g_{j,i,\g}:B(z_j,r)\rightarrow\C$ as 
$g_{j,i,\g}(z)=\frac{\varphi_i^\g(z)}{(\varphi_i^\g)'(z_j)}$. 
Then $g_{j,i,\g}$ is a univalent function on $B(z_j,r)$ such that 
$g_{j,i,\g}'(z_j)=1$. 
By Koebe's Distortion Theorem, there exists a constant $M_0\ge 1$ such that 
for all $1\leq j\leq J$, $i\in I$ and $\g\in D_d(\g_0,R)$, 
we have $M_0^{-1}\leq|g_{j,i,\g}'|\leq M_0$ on $B(z_j,r/2)$. 
Thus, there exists a constant $M_1\ge 1$ such that 
for all $i\in I$, $\g\in D_d(\g_0,R)$, and $a,b\in X$, 
$$M_1^{-1}\leq\left|\frac{(\varphi_i^\g)'(b)}{(\varphi_i^\g)'(a)}\right|\leq M_1.$$
Combining this with Theorem 4.1.2 in~\cite{gdms}, 
we obtain that there exists a constant $M_2>0$ such that 
for all $i\in I$, $\g\in D_d(\g_0,R)$ and 
$\om, \tau\in I^\infty$ such that $\om_1=\tau_1=i$, 
\begin{equation}
\label{sumiremeq3}
\aligned
|(\varphi_i^\g)'(\beta_{\g,\om})-(\varphi_i^\g)'(\beta_{\g,\tau})|
& \leq M_2|(\varphi_i^\g)'(\beta_{\g,\tau})| |\beta_{\g,\om}-\beta_{\g,\tau}|\\ 
& \leq M_2|(\varphi_i^\g)'(\beta_{\g,\tau})| e^{-\a|\om\wedge\tau|}\mbox{diam}(X).
\endaligned
\end{equation} 
From~(\ref{mar2}),~(\ref{sumiremeq3}) and the assumption $|\ka_\om(\g)|\leq B$, 
it follows that there exists a constant $M_3>0$ such that 
for all $\om, \tau\in I^\infty$ with $\om_1=\tau_1=i$ and for all $\g\in D_d(\g_0,R)$, 
$$|\ka_\om(\g)-\ka_\tau(\g)|\leq M_3 e^{-\a|\om\wedge\tau|}.$$
Applying Cauchy's formula to the well-defined branch  
$\log:B(1,1/2)\rightarrow\C$ with $\log 1=0$, which is bounded,    
we obtain that there exists a constant $M_4>0$ such that 
$$|\log\ka_\om(\g)-\log\ka_\tau(\g)|\leq M_4 e^{-\a|\om\wedge\tau|}$$
for all $\om, \tau\in I^\infty$ with $\om_1=\tau_1=i$ and for all $\g\in D_d(\g_0,R)$.
The claim holds with $C_0=\max\{2M,M_4\}$.

\

Since by Lemma~\ref{l1myu9.1} the operator $\Ra_d^{\g_0}:\Hol(D_b^d(\g_0,R))\to \Hol(D_b^{2d}(I(\g_0),R/4))$
is linear and bounded with norm bounded above by $4^d$, we get that for all 
$\^\g\in D_{2d}(I(\g_0),R/4)$, 
$$
\big|\Ra_d^{\g_0}(\log\ka_\om)(\^\g)-\Ra_d^{\g_0}(\log\ka_\tau)(\^\g)\big|
=\big|\Ra_d^{\g_0}(\log\ka_\om-\log\ka_\tau)(\^\g)\big|
\le 4^dC_{0} e^{-\a|\om\wedge\tau|}.
$$
Hence, $\Ra_d^{\g_0}(\log\ka)(\^\g)\in\Ha_\a$ and $v_\a(\Ra_d^{\g_0}(\log\ka)(\^\g))
\le 4^dC_{0}$ for all $\^\g\in D_{2d}(I(\g_0),R/4)$. Combining this with (\ref{1rs259}) 
we deduce that $\|\Ra_d^{\g_0}(\log\ka)(\^\g)\|_\a
\le 4^d(M+C_{0})$. Thus Condition (a) of Lemma~\ref{l2rs253} is satisfied. The proof
is complete. \endpf

\

We are finally in a position to prove our theorem.

\

{\bf Proof of Theorem~\ref{prealanalconj}.} Let $\g_0$, $R_0$, $s$, $\a$, and $R$ 
be as in Lemma~\ref{l1rs259}. Let further $\^\Ga=\C^{\,2d}$, 
$G=D_{2d}(I(\g_0),R/4)$, $\^\g_0=I(\g_0)$, $\Phi^{\^\g_0}=\Phi^{\g_0}$, and  $g_{\^\g}(\om)=t\Ra_d^{\g_0}(\log\kappa_\om)(\^\g)$, as in the proof of Lemma~\ref{l1rs259}. 
Recall that according to that lemma, the operator $\L_{\^\g}:\Ha_\a\to\Ha_\a$ generated 
by the potential $u_{\^\g}(\om):=t\Ra_d^{\g_0}(\log\ka_\om)(\^\g)
+t\log|(\varphi_{\om_1}^{\g_0})'(\pi_{\g_0}(\sg(\om)))|$ depends holomorphically on 
$\^\g\in D_{2d}(I(\g_0),R/4)$. Note that, because of the 
definition of the function $\ka$, we have
\begin{equation}\label{1rs041408}
\exp\(t\Ra_d^{\g_0}(\log\ka_\om)(\g)+t\log|(\varphi_{\om_1}^{\g_0})'(\pi_{\g_0}(\sg(\om)))|\)
=|(\varphi_{\om_1}^{\g})'(\pi_{\g}(\sg(\om)))|^t
\end{equation}
for all $\g\in D_{d}(\g_0,R/4)$, where $D_{d}(\g_0,R/4)$ is identified with $D_{2d}(I(\g_0),R/4)\cap\R^{2d}$ as 
$\C^d$ is embedded in $\C^{\,2d}$ by the formula $(x_1+iy_1,x_2+iy_2,\ld,x_d+iy_d)
\mapsto (x_1,y_1,x_2,y_2,\ld,x_d,y_d)$. Thus, in view of (\ref{1rs041408}) and 
Theorem~2.4.6 and Corollary~2.7.5 from \cite{gdms}, $e^{P(\g_0,t)}$ is a simple isolated eigenvalue 
of the operator $\L_{I(\g_0)}:\Ha_\a\to\Ha_\a$. Hence, in view of Lemmas~\ref{l2rs253}
and~\ref{l1rs259}, which ensure the analyticity of the function $\^\g\mapsto 
\L_{\^\g}$ from $D_{2d}(I(\g_0),R/4)$ to $L(\Ha_\a)$,
Kato-Rellich Perturbation Theorem~(see~\cite{katorell}, 
Theorem 3.16 on page 212 and pages 368-369) yields $r_1\in(0,R/4]$ and
a holomorphic function $\xi:D_{2d}(I(\g_0),r_1)\to\C$ such that $\xi(I(\g_0))
=e^{P(\g_0,t)}$ and such that for each $\^\g\in D_{2d}(I(\g_0),r_1)$, 
it turns out that $\xi(\^\g)$ is a simple isolated eigenvalue of $\L_{\^\g}:
\Ha_\a\to\Ha_\a$ with the remainder of the spectrum uniformly separated from
$\xi(\^\g)$. In particular, there exists $r_2\in(0,r_1]$ and $\eta>0$ such that
\begin{equation}\label{1120905p90}
\sg(\L_{\^\g})\cap B(e^{P(\g_0,t)},\eta)=\{\xi(\^\g)\}
\end{equation}
for all $\^\g\in D_{2d}(I(\g_0),r_2)$. Since $e^{P(\g,t)}$ is the spectral radius 
of the operator $\L_{\g}$ for all $\g\in D_d(\g_0,r_2)$ 
(see Theorem~2.4.6 and Corollary~2.7.5 in~\cite{gdms}), in view 
of the semi-continuity of the spectral set function (see Theorem~10.20 
in~\cite{rudin}), taking $r_2>0$ small enough, we also have that 
the spectral radius of $\L_\g$ is in $[0,e^{P(\g_0,t)}+\eta)$, and 
along with~(\ref{1120905p90}), this implies that $e^{P(\g,t)}=\xi(\g)$ for all
$\g\in D_d(\g_0,r_2)$. Consequently, the function $\g\mapsto P(\g,t)$, 
$\g\in D_d(\g_0,r_2)$, is real-analytic. The proof of Theorem~\ref{prealanalconj}
is complete. \endpf 

\
                   
\begin{rem}
Condition~(1) in Theorem~\ref{prealanalconj} is a local 
version of condition~(c) in~\cite{RU}. 
\end{rem}

Those analytic families that depend continuously 
on the parameter $\g$ when $\HIFS(X)$ is equipped with the 
$\lambda$-topology satisfy condition~(1). 

\begin{lem}\label{lambdacontilem}
Let $\{\Phi^\g\}_{\g\in\Ga}$
be an analytic family such that 
$\g\mapsto\Phi^\g\in\HIFS(X)$ is 
continuous with respect to the $\lambda$-topology.
Then condition~(1) in Theorem~\ref{prealanalconj} holds. 
\end{lem}

{\sl Proof.} This lemma is a consequence of Lemma~\ref{l:pigomegaanal}  
and the observation following Corollary~6.2 in~\cite{RU}. 
\endpf 


\

Remark~\ref{rnc}, Theorems~\ref{psubharm} and~\ref{prealanalconj}, and Lemma~\ref{lambdacontilem}
justify the assumption of continuous dependence on $\g$ of the families studied in the sequel.
 

\section{A Classification Theorem}\label{ct}

We now give a complete classification of analytic families which are 
continuous with respect to the $\lambda$-topology.

\bdfn
Let $\{\Phi^\g\}_{\g\in\Ga}$ 
be an analytic family in $\HIFS(X)$. 
We define the following subsets of parameters:
\begin{itemize}
\item 
$\RE_\Ga(X):=\{\g\in\Ga: 
\Phi^\g\in\RE(X)\}$;
\item 
$\SR_\Ga(X):=\{\g\in\Ga: 
\Phi^\g\in\SR(X)\}$; 
\item 
$\CFR_\Ga(X):=\{\g\in\Ga:\Phi^\g\in\CFR(X)\}$;
\item 
$\FSR_\Ga(X):=\{\g\in\Ga:\Phi^\g\in\SR(X)\sms\CFR(X)\}$;
\item 
$\CR_\Ga(X):=\{\g\in\Ga:\Phi^\g\in\CR(X)\}$;
\item 
$\IR_\Ga(X):=\{\g\in\Ga:\Phi^\g\in\IR(X)\}$; 
\item 
$\RAH(\Ga):=\{\g_0\in\Ga: 
\g\mapsto h_{\Phi^\g}\mbox{ is real-analytic in a neighborhood of } 
\g_0\}$;
and   
\item 
$\NPHH(\Ga):=\{\g_0\in\Ga:\g\mapsto 
h_{\Phi^\g}\mbox{ is not pluriharmonic in any neighborhood of } 
\g_0\}$.
\end{itemize}
\edfn

\

\bthm\label{classthm}
Let $\{\Phi^\g\}_{\g\in\Ga}$ 
be an analytic family in $\HIFS(X)$ such that 
$\g\mapsto\Phi^\g\in\CIFS(X)$ is continuous with respect to 
the $\lambda$-topology. Let $\theta$ be the constant such that 
$\theta_{\Phi^\g}\equiv\theta$ on $\Ga$ (this constant exists 
according to Lemma~5.4 in~\cite{RU}). 
Then exactly one of the following statements holds. 
\begin{itemize}
\item[(I)]
$\Ga=\CFR_\Ga(X)=\RAH(\Ga)$, and the function $\g\mapsto 
h_{\Phi^\g}$, $\g\in\Ga$, is identically equal to 
a constant larger than $\theta$.  
\item[(II)]
$\Ga=\CFR_\Ga(X)=\RAH(\Ga)=\NPHH(\Ga)$.   
\item[(III)]
$\Ga=\FSR_\Ga(X)=\RAH(\Ga)$, 
and the function $\g\mapsto h_{\Phi^\g}$, $\g\in\Ga$, 
is identically equal to a constant larger than $\theta$.
\item[(IV)]
$\Ga=\FSR_\Ga(X)=\RAH(\Ga)=
\NPHH(\Ga)$. 
\item[(V)] 
$\Ga=\FSR_\Ga(X)\cup\CR_\Ga(X) 
=\NPHH(\Ga)$, with $\FSR_\Ga(X)\neq 
\emptyset$ and $\CR_\Ga(X)\neq\emptyset$. 

Furthermore, 
$\CR_\Ga(X)=\partial(\FSR_\Ga(X))$. 

Also, $\FSR_\Ga(X)\subset\RAH(\Ga)$. 
\item[(VI)] 
$\Ga=\CR_\Ga(X)=\RAH(\Ga)$, and the function 
$\g\mapsto h_{\Phi^\g}$, $\g\in\Ga$,  is 
identically equal to the constant $\theta$. 
\item[(VII)]
$\FSR_\Ga(X)\neq 
\emptyset$, $\CR_\Ga(X)\neq\emptyset$ 
and $\IR_\Ga(X)\neq\emptyset$, while $\CFR_\Ga(X)=\emptyset$. 

Moreover, $\emptyset\neq\partial(\IR_\Ga(X))
\subset\CR_\Ga(X)=\partial(\FSR_\Ga(X))$. 

Moreover,
$\FSR_\Ga(X)\cup\IR_\Ga(X)
\subset\RAH(\Ga)$. 

Moreover, 
$\emptyset\neq\partial(\IR_\Ga(X))\subset\Ga\sms 
\RAH(\Ga)$.

Furthermore, 
$\NPHH(\Ga)=\FSR_\Ga(X)\cup\CR_\Ga(X)$.  
\item[(VIII)]
$\Ga=\IR_\Ga(X)=\RAH(\Ga)$, and the function  
$\g\mapsto h_{\Phi^\g}$, $\g\in\Ga$, 
is identically equal to the constant $\theta$. 
\end{itemize}
\ethm

{\sl Proof.} 
{\bf Case (1):} 
Suppose that $\CFR_\Ga(X)\neq\emptyset$. Since 
$\CFR_\Ga(X)$ is clopen by Lemma~5.9 in~\cite{RU}
and $\Ga$ is connected, 
we deduce that $\Ga=\CFR_\Ga(X)$. Now, if 
the function $\g\mapsto h_{\Phi^\g}$, $\g\in\Ga$, is not 
constant, then combining Theorem~6.1 and Corollary~6.4 in~\cite{RU} 
we obtain $\Ga=\RAH(\Ga)=\NPHH(\Ga)$. This corresponds to statement (II). 

On the other hand, if the function $\g\mapsto h_{\Phi^\g}$, $\g\in\Ga$,
is identically equal to a constant $\tau$, then by Theorem~2.4 in~\cite{RU} 
we have $\tau>\theta$. This is statement (I). 

{\bf Case (2):} Suppose next that 
$\Ga=\FSR_\Ga(X)$. 
If $\g\mapsto h_{\Phi^\g}$, $\g\in\Ga$, is not 
constant, then combining Theorem~6.1 and Corollary~6.4 in~\cite{RU} 
we get $\Ga=\RAH(\Ga)=\NPHH(\Ga)$. This is statement (IV).  

If, however, the function $\g\mapsto h_{\Phi^\g}$, $\g\in\Ga$,
is identically equal to a constant $\tau$, then by Theorem~2.4 in~\cite{RU} we have
$\tau>\theta$. This is statement (III). 

{\bf Case (3):} 
Suppose that $\Ga=\FSR_\Ga(X)\cup\CR_\Ga(X)$,  
that $\FSR_\Ga(X)\neq\emptyset$, 
and $\CR_\Ga(X)\neq\emptyset$. 
Then, by Lemma~5.9 in~\cite{RU} we know that $\FSR_\Ga(X)$ 
is a proper non-empty open subset of $\Ga$. 
Moreover, by Theorem~6.1 in~\cite{RU} we have 
$\FSR_\Ga(X)\subset\RAH(\Ga)$. 
Now, if there exists a parameter $\g_0\in\FSR_\Ga(X)$ 
such that $\g\mapsto h_{\Phi^\g}$ is pluriharmonic around 
$\g_0$, then by Corollary~6.4 and Theorem~2.4 in~\cite{RU}, we deduce that 
the function $\g\mapsto h_{\Phi^\g}$ is identically equal to 
a constant $\tau>\theta$ in the connected component $U_0$ of 
$\FSR_\Ga(X)$ containing $\g_0$. 
Since $\CR_\Ga(X)\neq\emptyset$ by assumption and since 
$\Ga $ is connected, the boundary $\partial U_0$ of $U_0$ 
in $\Ga$ is not empty. 
Moreover, $\partial U_0\subset\CR_\Ga(X)$.  
Then, according to Theorem~5.10 in~\cite{RU}, the function 
$\g\mapsto h_{\Phi^\g}$ is continuous and therefore 
$h_{\Phi^\g}=\tau$ for every $\g\in\partial U_0$.
However, using Theorem~2.3 in~\cite{RU},
we know that $h_{\Phi^\g}=\theta$ for every $\g\in\CR_\Ga(X)$.
This is a contradiction.  Consequently, 
$\overline{\FSR_\Ga(X)}
\subset\NPHH(\Ga)$. Furthermore, Theorem~\ref{prealanalconj} above  
and Theorems~2.3 and~2.4 in~\cite{RU} imply that $\CR_\Ga(X)=
\partial(\FSR_\Ga(X))$.   
Thus, we obtain $\Ga=
\overline{\FSR_\Ga(X)}=\NPHH(\Ga)$. 
This means that, in this case, statement (V) holds. 

{\bf Case (4):} 
Suppose next that $\Ga=\CR_\Ga(X)$. Then, 
by Theorem~2.3 in~\cite{RU}, statement (VI) holds. 

{\bf Case (5):} 
Suppose now that 
$\Ga=\IR_\Ga(X)$. Then, 
by Theorem~2.3 in~\cite{RU}, statement (VIII) holds. 

{\bf Case (6):} 
Suppose next that 
$\CR_\Ga(X)\neq\emptyset$ and 
$\IR_\Ga(X)\neq\emptyset$. 
Then $\g\mapsto P(\g,\theta)$ is not constant on 
$\Ga$. Theorems~\ref{prealanalconj} and~\ref{psubharm} imply that 
the function $\g\mapsto P(\g,\theta)$ is  
a non-constant real-analytic function that satisfies the Maximum Principle. 
It follows that for any parameter $\g_0\in\CR_\Ga(X)$, 
there exists a parameter $\g$, arbitrarily close to 
$\g_0$, such that $P(\g,\theta)>0$. Notice that
$P(\g,\theta)>0$ implies that $\g\in\SR_\Ga(X)$. Hence,   
$\CR_\Ga(X)\subset\partial(\FSR_\Ga(X))$. 
On the other hand, the continuity of the function 
$\g\mapsto P(\g,\theta)$ implies that
$\partial(\FSR_\Ga(X))\subset\CR_\Ga(X)$.
Thus, $\partial(\FSR_\Ga(X))=\CR_\Ga(X)$.
The continuity of the function $\g\mapsto P(\g,\theta)$ also 
implies that $\partial(\IR_\Ga(X))\subset\CR_\Ga(X)$, while the 
connectedness of $\Ga$ ensures that $\partial(\IR_\Ga(X))\neq\emptyset$.  

Moreover, using Theorems~5.10 and~6.1 in~\cite{RU}, 
we deduce that $\FSR_\Ga(X)
\cup\IR_\Ga(X)\subset\RAH(\Ga)$. 

In order to show that $\partial(\IR_\Ga(X))\subset 
\Ga\sms\RAH(\Ga)$, 
let $\g_0\in\partial(\IR_\Ga(X))$. 
Then, $\g_0\in\CR_\Ga(X)=
\partial(\FSR_\Ga(X))$. 
If $\g_0\in\RAH(\Ga)$, then 
the function $\g\mapsto h_{\Phi^\g}$ is identically 
equal to the constant $\theta$ in a neighborhood of $\g_0$. 
However, this contradicts 
Theorem~2.4 in~\cite{RU}. Hence, we obtain that 
 $\partial(\IR_\Ga(X))\subset 
\Ga\sms\RAH(\Ga)$. 

By the same argument as in case (3), we obtain 
$\NPHH(\Ga)=\FSR_\Ga(X)
\cup\CR_\Ga(X)$. 

Thus, in this case, statement 
(VII) holds. 

{\bf Case (7):} Finally, suppose that 
$\FSR_\Ga(X)\neq\emptyset$ 
and $\IR_\Ga(X)\neq\emptyset$. 
Then, by Theorem~5.7 in~\cite{RU}, we have $\CR_\Ga(X)\neq\emptyset$. 
Therefore, by case (6), statement (VII) holds. 

We are done.     
\endpf

\

This classification theorem further highlights some differences between  
analytic families and the entire space $\CIFS(X)$. 
Whereas we know that $\partial(\IR(X))=\CR(X)$ while $\partial(\SR(X))\subset\CR(X)$ 
in $\CIFS(X)$, there exist analytic families of types~(V) and~(VII) (see Example~\ref{example}, Theorem~\ref{t:hsptdense}, 
and Proposition~\ref{p:everytype}).

\bcor\label{classcor}
Under the assumptions of Theorem~\ref{classthm}, 
the set $\Ga\sms\RAH(\Ga)$ is included in a proper 
real-analytic subvariety of $\Ga$. In particular, 
$\RAH(\Ga)$ is open and dense in $\Ga$.
\ecor

{\sl Proof.} 
If the family $\{\Phi^\g\}_{\g\in\Ga}$ is 
of type (V) or (VII), then 
$\Ga\sms\RAH(\Ga)$ is included in 
the set $\CR_\Ga(X)=\{\g\in\Ga\mid P(\g,\theta)=0\}$, 
which is a proper real-analytic subset of $\Ga$ due to
Theorem~\ref{prealanalconj}.
\endpf



\

Note that every type of family described in Theorem~\ref{classthm} 
exists (see Proposition~\ref{p:everytype}). It is indeed possible to construct families of every type 
and in the following example we present a family of type (V) in $\HIFS(\overline{\D})$, where
$\D=B(0,1)\sbt\C$ is the open unit disk and $\overline{\D}=\overline{B(0,1)}\sbt\C$ is the closed unit disk. 

\bex\label{example}
Let $U$ be a non-empty open subset of $\D$ and 
let $z_1,z_2,z_3\in U$ be mutually distinct points. 
Let $\{\alpha_i(z):=a_i(z-z_i)+z_i\}_{i=1}^3$ be  
such that $0<a_1=a_2=a<1/2$, $0<a_3<1$, 
$\bigcup_{i=1}^3\alpha_i(\overline{\D})\subset U$, 
and $\alpha_i(\overline{\D})\cap\alpha_j(\overline{\D})=\emptyset$ whenever 
$i\neq j$. 
Let $\{c_i\}_{i\geq 3}$ be a sequence of positive numbers such that 
$$\sum_{i=3}^\infty c_i^t=
\left\{\begin{array}{ll}
\infty & \mbox{ if } t<1 \\ 
\frac{1-2a}{a_3}      & \mbox{ if } t=1 \\
<\infty& \mbox{ if } t>1. 
\end{array}\right.
$$ 
Let $\Ga=B(0,r)\sbt\C$ be a small disk around $0$ and 
$X=\overline{\D}$. 
For each $\g\in\Ga$, let 
$$\Phi^\g:= 
\begin{cases}
\varphi_1^\g(z):=4a(\frac{1}{2}+\g)^{2}(z-z_1)+z_1,\\ 
\varphi_2^\g(z):=4a(\frac{1}{2}-\g)^{2}(z-z_2)+z_2,\\ 
\varphi_i^\g(z):=a_3(\psi_i(z)-z_3)+z_3, \hspace{1cm} i\geq 3,
\end{cases}
$$
where $\{\psi_i\}_{i\geq 3}$ is a family of 
similitudes  such that 
$\psi_i(\D)\subset\D$, 
$\psi_i(\D)\cap\psi_j(\D)=\emptyset$ whenever $i\neq j$, 
and $\psi_i'(z)\equiv c_i$. (Note that 
by Lemma 3.2(i) in~\cite{RSU}, 
such $\{c_i\}_{i\geq 3}$ and $\{\psi _{i}\}_{i\geq 3}$ exist.)  
We claim that the family $\{\Phi^\g\}_{\g\in\Ga}$ is of type (V) 
and $\varphi_i^\g(\overline{\D})\subset U$ for all $\g\in\Ga$ and all $i\in I$. 
\eex 

Indeed, if we take $r$ so small, $\{\Phi^\g\}_{\g\in\Ga}$ is an analytic 
family in $\HIFS(\overline{\D})$ and 
$\varphi_i^\g(\overline{\D})\subset U$ for all $\g\in\Ga$ and all $i\in I$. 
Moreover, $\g\mapsto\Phi^\g$ is continuous with respect to the $\lambda$-topology. 

Furthermore, 
$$P(\g,t)=\log\left((4a)^t\Bigl|\frac{1}{2}+\g\Bigr|^{2t}+
(4a)^t\Bigl|\frac{1}{2}-\g\Bigr|^{2t}
+a_3^t\sum_{i=3}^{\infty}c_i^t\right).  
$$
Therefore $\th_{\Phi^\g}=1$ for every $\g$ and 
\begin{eqnarray*}
P(\g,\th_{\Phi^\g})=
P(\g,1)&=&\log\left(4a\Bigl|\frac{1}{2}+
\g\Bigr|^2+4a\Bigl|\frac{1}{2}-\g\Bigr|^2+a_3\sum\limits_{i=3}^{\infty }c_i\right) \\
&=&\log\left(4a\Bigl(\frac{1}{2}+2|\g|^2\Bigr)+a_3\cdot\frac{1-2a}{a_3}\right) \\ 
&=&\log\left(1+8a|\g|^2\right).
\end{eqnarray*}
This latter function of $\g$ takes its minimum value $0$ at $\g=0$ 
and is positive for each small $\g\neq 0$. 
Thus, $\{\Phi^\g\}_{\g\in\Ga}$ is of 
type (V).
\endpf 

\section{Consequences of the Classification Theorem}\label{cct}

Propositions~5.13--5.16 in~\cite{RSU} assert that
$\partial(\SR(X))=\partial(\FSR(X))\sbt\CR(X)$, 
despite that $\CR(X)\neq\partial(\SR(X))$. Moreover, Proposition~5.17 
of that same paper provides a condition under which 
a critically regular SIFS is in the boundary of the subset 
of strongly regular CIFSs. Using our classification theorem, 
we partially generalize this result.   
In this section, let $X$ be a non-empty compact connected subset of $\C $ 
such that $X=\overline{\mbox{Int}(X)}$ and 
the condition (iii) in section~\ref{secprelifs} is satisfied.
\bthm\label{crxthm}
Suppose that $X\subset\C$ is star shaped with $x_0\in\Int(X)$ for center. 
Let $\Phi=\{\varphi_{i}\}\in\CR(X)\cap\HIFS(X)$ be a system   
such that there exists a finite set $F\sbt\N$ for which 
$$\bigcup_{i\in F}\varphi_{i}(X)\subset\Int(X)\sms\overline{\bigcup_{j\not\in F}\varphi_{j}(X)}
\neq \emptyset.$$ 
Then $\Phi\in\partial(\SR(X))=\partial(\FSR(X))$. More precisely, 
there exists an open, connected neighborhood $\Ga$ of $\{s\in\R: 
0<s\leq 1\}$ in $\C$ such that, 
setting 
$$\Phi^\g=\left\{\begin{array}{ll}
                            \varphi_i, & i\not\in F, \\
                            \varphi_i\circ\psi_\g, & i\in F
                      \end{array}
               \right\}
$$                            
for each $\g\in\Ga$, where $\psi_{\g}(x)=\g(x-x_{0})+x_{0}$,  
the analytic family $\{\Phi^\g\}_{\g\in\Ga}$ 
in $\HIFS(X)$ satisfies all of the following: 
\begin{itemize}
\item[(1)]
$\Phi^{1}=\Phi$; 
\item[(2)] $\g\in\Ga\mapsto\Phi^\g\in \HIFS(X)$ is continuous 
with respect to the $\lambda$-topology; 
\item[(3)] 
$\{\Phi^\g\}_{\g\in\Ga}$ is of type 
(VII); 
\item[(4)]
$1\in\CR_\Ga(X)=\partial(\FSR_\Ga(X))$;  
\item[(5)] 
$\{s\in\R:0<s<1\}\cap\IR_\Ga(X)\neq\emptyset$ 
and $\{s\in\R:0<s\leq 1\}\cap\partial(\IR_\Ga(X))\neq\emptyset$; and  
\item[(6)] 
for any $\g_{0}\in\partial(\IR_\Ga(X))$, 
we have $\g_{0}\in\partial(\IR_\Ga(X))\subset 
\CR_\Ga(X)=\partial(\FSR_\Ga(X))$ 
and $\g\mapsto h_{\Phi^\g}$ is not 
real-analytic in any neighborhood of $\g_{0}$ in $\Ga$. 
\end{itemize}  
\ethm

{\sl Proof.}  
Using the local 
bounded distortion property (see Corollary~4.1.4 in~\cite{gdms}), 
the fact that all generators $\varphi_i$, $i\in\N$, are
contractions with a common ratio, and the Mean Value Inequality,
we can construct a decreasing family $\cal F$ of open, connected 
neighborhoods of $X$ whose intersection is $X$ and with respect 
to each of which $\Phi$ is in $\HIFS(X)$. (In fact, $X_\e:=B(X,\e)$ 
is such a neighborhood for every $\e>0$
small enough). This family $\cal F$ therefore has a member $U\supset X$
such that 
$$\bigcup_{i\in F}\varphi_{i}(\overline{U})
\subset\Int(X)\sms\overline{\bigcup_{j\not\in F}\varphi_{j}(X)}
\neq \emptyset.$$ 
Choose $\delta>0$ such that $\overline{X_\delta}\subset U$.
Thereafter choose an open neighborhood $\Ga$ of 
$\{s\in\R:0<s\leq 1\}$ in $\C$ such that 
\[ \bigcup_{\g\in\Ga}\psi_\g(\overline{X_\delta})
\subset U. \]
Finally, take for $V$ any neighborhood in the family $\cal F$ such that
$V\subset X_\delta$.
Then $\varphi_i(V)\subset V$ for all $i\in\N$ and 
\[ \varphi_k\circ\psi_\g(V)
   \subset\varphi_k\circ\psi_\g(\overline{X_\delta})
   \subset\varphi_k(U)
   \subset\Int(X)\sms\overline{\cup_{j\not\in F}\varphi_{j}(X)} 
\]  
for every $k\in F$ and every $\g\in\Ga$. With this neighborhood $V$, 
we have that $\{\Phi^\g\}_{\g\in\Ga}$ 
is an analytic family in $\HIFS(X)$. Moreover, 
$\g\mapsto\Phi^{\g}\in \HIFS(X)$ is 
clearly continuous with respect to the $\lambda$-topology.  
Furthermore, if $s>0$ is small enough, then by the argument given 
in the proof of Proposition~5.11 in~\cite{RSU}, 
we obtain that $s\in\IR_\Ga(X)$. 
Thus, Theorem~\ref{classthm} implies that 
the analytic family $\{\Phi^\g\}_{\g\in\Ga}$ 
is of type (VII). In particular, $1\in\CR_\Ga(X)=
\partial(\FSR_\Ga(X))$. 

Now, let 
\[ S:=\sup\Bigl\{s\in\R:0<s<1,\,s\in\IR_\Ga(X)\Bigr\}. \]
Then $S\in\partial(\IR_\Ga(X))$.  
Thus, statement (5) holds. 

Finally, by Theorem~\ref{classthm}, 
statement (6) holds. 
\endpf

\

By a similar argument, we get the following, more general result. 

\bthm\label{crxthm2}
Let $X\subset\C$ and 
$\Phi=\{\varphi_{i}\}\in\CR(X)\cap\HIFS(X)$ be such that 
there exist a finite set $F\sbt\N$ 
and a simply connected subdomain $W$ of $\Int(X)$ such that 
\[ \bigcup_{i\in F}\varphi_{i}(X)\subset W \mbox{ and } 
W\bigcap\overline{\bigcup_{j\not\in F}\varphi_{j}(X)}=\emptyset. \] 
Then $\Phi\in\partial(\SR(X))=\partial(\FSR(X))$. More precisely, 
there exists an analytic family 
$\{\Phi^\g\}_{\g\in\Ga}$ 
in $\HIFS(X)$ with an open, connected neighborhood 
$\Ga$ of $\{s\in\R:0<s\leq 1\}$ in $\C$ such that 
\begin{itemize}
\item[(1)]
$\Phi^{1}=\Phi$; 
\item[(2)]
$\g\in\Ga\mapsto\Phi^\g\in\HIFS(X)$ is continuous 
with respect to the $\lambda$-topology; 
\item[(3)]
$\{\Phi^\g\}_{\g\in\Ga}$ is of type 
(VII); 
\item[(4)]
$1\in\partial(\FSR_\Ga(X))$;
\item[(5)]
$\{s\in\R:0<s<1\}\cap\IR_\Ga(X)\neq\emptyset$ 
and $\{s\in\R:0<s\leq 1\}\cap\partial(\IR_\Ga(X))
\neq\emptyset$; and  
\item[(6)]
for any $\g_{0}\in 
\partial(\IR_\Ga(X))$, 
we have $\g_{0}\in\partial(\IR_\Ga(X))\subset 
\CR_\Ga(X)=\partial(\FSR_\Ga(X))$ 
and $\g\mapsto h_{\Phi^\g}$ is not 
real-analytic in any neighborhood of $\g_{0}$ in $\Ga$. 
\end{itemize}
\ethm

{\sl Proof.} Let $A$ be a simply connected subdomain of $W$ 
such that $\cup_{i\in F}\varphi_{i}(X)\subset 
A\subset\overline{A}\subset W$. 
Let also $a\in A$. Taking an appropriate $A$, 
the Riemann Mapping Theorem implies that there exists an analytic 
family $\{\psi_{\g}\}_{\g\in\Ga}$ 
of injective holomorphic maps from $A$ to $W$ 
with an open, connected neighborhood $\Ga$ of 
$\{s\in\R:0<s\leq 1\}$ in $\C$ 
such that $\psi_{1}(z)\equiv z$ and  
$\psi_{s}(z)\rightarrow a$ as $|s|\rightarrow 0$ 
uniformly on $A$. 

For each $\g\in\Ga$, let 
$$\Phi^\g=\left\{\begin{array}{ll}
                            \varphi_i, & i\not\in F, \\
                            \psi_\g\circ\varphi_i, & i\in F
                      \end{array}
               \right\}
$$                            
Then 
$\{\Phi^{\g}\}_{\g\in\Ga}$ is an 
analytic family in $\HIFS(X)$. 
Note that a neighborhood common to all $\Phi^\g$'s
can be found as follows. As explained in the proof of Theorem~\ref{crxthm}, 
we can construct a decreasing family $\cal F$ of open, connected 
neighborhoods of $X$ whose intersection is $X$ and with respect 
to each of which $\Phi$ is in $\HIFS(X)$. This family 
therefore has a member $V\supset X$ such that 
$\cup_{i\in F}\varphi_i(V)\subset A$.
Then $\varphi_i(V)\subset V$ for all $i\in\N$ and 
\[ \psi_\g\circ\varphi_i(V)
   \subset\psi_\g(A)
   \subset W
   \subset \Int(X)
   \subset V\]
for every $i\in F$. Moreover, 
$\g\in\Ga\mapsto\Phi^\g\in\HIFS(X)$ is 
continuous with respect to the $\lambda$-topology.  
Furthermore, by the argument given in the proof of Proposition~5.11
in~\cite{RSU}, every $\g\in\Ga$ with small enough modulus 
is in $\IR_\Ga(X)$. 
Hence, the family $\{\Phi^\g\}_{\g\in\Ga}$ 
is of type (VII).
In particular, $1\in\CR_\Ga(X)=\partial(\FSR_\Ga(X))$. 
The remaining statements follow from the same arguments 
as those given in the proof of Theorem~\ref{crxthm}.
\endpf

\

We now give a better description of the set 
$\FSR(X)\cap\HIFS(X)$. 

\bprop\label{propo}
Suppose that $X$ is star shaped with $x_0\in\Int(X)$ for center. 
%
Let $\Phi=\{\varphi_{i}\}\in\FSR(X)\cap\HIFS(X)$. 
Then there exists an open, connected neighborhood 
$\Ga$ of $\{s\in\R:0<s\le 1\}$ in $\C$ 
such that, setting $\Phi^\g=
\{\varphi_{i}\circ\psi_{\g}\}_{i\in I}$ 
for each $\g\in\Ga$, where 
$\psi_{\g}(x)=\g(x-x_{0})+x_{0}$, 
the analytic family 
$\{\Phi^\g\}_{\g\in\Ga}$ 
in $\HIFS(X)$ satisfies 
\begin{itemize}
\item[(1)]
$\Phi^{1}=\Phi$; 
\item[(2)] $\Phi^\g\rightarrow\Phi$ 
as $\g\nearrow 1$ in $(0,1)$ 
with respect to the $\lambda$-topology; 
\item[(3)]
$\g\in\Ga\mapsto\Phi^\g\in\HIFS(X)$ is 
continuous with respect to the $\lambda$-topology; 
\item[(4)] 
$\{\Phi^\g\}_{\g\in\Ga}$ is of type 
(VII); 
\item[(5)] 
$\{s\in\R:0<s<1\}\cap\IR_\Ga(X)\neq\emptyset$ 
and $\{s\in\R:0<s<1\}\cap\partial(\IR_\Ga(X))
\neq\emptyset$; and  
\item[(6)] 
for any $\g_{0}\in\partial(\IR_\Ga(X))$, 
we have $\g_{0}\in\partial(\IR_\Ga(X))\subset 
\CR_\Ga(X)=\partial(\FSR_\Ga(X))$ 
and $\g\mapsto h_{\Phi^\g}$ is not 
real-analytic in any neighborhood of $\g_{0}$ in $\Ga$. 
\end{itemize}
\eprop

{\sl Proof.} 
It is easy to see that the first three statements hold. 
Moreover, if $\g\in\Ga$ and $|\g|$ is small enough, 
then $\g\in\IR_\Ga(X)$. 
Since $\g\in\SR_\Ga(X)$ if $\g$ is close enough to 
$1$, statement (4) follows from Theorem~\ref{classthm}. 
By the same argument as that in the proof of 
Theorem~\ref{crxthm}, we can show that statements (5) and (6) 
hold.   
\endpf

%

\

We now take a brief look at a subset of $\CR(X)\cap\HIFS(X)$.

\bdfn
Let $X\subset\C$. 
Let $\HSPT(X)$ be the set of $\Phi\in\CR(X)\cap\HIFS(X)$ 
such that there exists an analytic family 
$\{\Phi^\g\}_{\g\in\Ga}$ in $\HIFS(X)$ 
with $\Ga=\{z\in\C:|z|<1\}$ 
satisfying all of the following:
\begin{itemize}
\item $\Phi^{0}=\Phi$; 
\item $\g\mapsto\Phi^\g\in\HIFS(X)$ is 
continuous with respect to the $\lambda$-topology; 
\item $\{\Phi^\g\}_{\g\in\Ga}$ 
is of type (VII); 
\item $0\in\partial(\IR_\Ga(X))
\subset \CR_\Ga(X)=
\partial(\FSR_\Ga(X))$; 
and 
\item $\g\mapsto h_{\Phi^\g}$ is not real-analytic 
in any neighborhood of $0$.  
\end{itemize}
\edfn

``HSPT'' stands for ``{\bf H}olomorphic {\bf S}imultaneous {\bf P}hase {\bf T}ransition''. 

\bthm
\label{t:hsptdense}
Let $X\subset\C$. 
Then $\HSPT(X)$ is a dense subset of $\HIFS(X)$ when this latter is 
endowed with the metric of pointwise convergence.
\ethm

{\sl Proof.} 
This proposition follows immediately from combining 
the argument in the proof of Lemma~4.10 in~\cite{RSU} 
(taking the similitude $S^{n}$ so that 
$S^{n}(X)$ is included in a small ball in $\varphi_{n}(X)$)
with Theorem~\ref{classthm}. 
\endpf

\ 

Finally, we show that every type of analytic family described in Theorem~\ref{classthm} exists. 
\bprop
\label{p:everytype}
Let $X\subset\C$. Then every type of analytic family in $\HIFS(X)$ described in 
Theorem~\ref{classthm} exists. 
\eprop

{\sl Proof.}
By Lemma 3.2(i) in~\cite{RSU}, there exist families of type 
(I), (III), (VI) and (VIII).
By Theorem~\ref{t:hsptdense}, there exists a family of type (VII).  
Taking an appropriate subfamily of type (VII), we obtain a family of type (IV).  
In order to construct a family of type (II), 
let $\Phi=\{\varphi_i\}_{i\in I}\in\SIFS(X)\cap\CFR(X)$.  
Let $z_0\in\mbox{Int}(X)$ and let 
$\alpha_\g(z)=\g(z-z_0)+z_0$ for each $\g\in\C$. 
Let $r>0$ be a small number such that 
for each $\g\in B(0,r)$, we have
$\alpha_\g(X)\subset\mbox{Int}(X)$. 
Let $\Ga=B(0,r)\backslash\{0\}$ and 
let $\Psi^\g=\{\alpha_\g\circ\varphi_i\}_{i\in I}$ for each 
$\g\in\Ga$. 
Then $\{\Psi^\g\}_{\g\in\Ga}$ is an analytic family in 
$\HIFS(X)$ such that $\g\mapsto\Psi^\g$ is continuous with respect to 
the $\lambda$-topology. Moreover, 
$P(\g,t)=P_\Phi(t)+t\log|\g|$ for all $\g\in\Ga$ and all $t\geq 0$. 
Hence, $\{\Psi^\g\}_{\g\in\Ga}$ is of type (II).  

We now construct a family of type (V). 
Let $S$ be a similarity of $\C$ such that 
$S(X)\subset\D$. 
Let $U:=\mbox{Int}(S(X))$ and let $z_1,z_2,z_3\in U$ be 
mutually distinct points. 
Let $\{\Phi^\g\}_{\g\in\Ga}$ 
be an analytic family of type (V) in $\HIFS(\overline{\D})$ such that 
$\varphi_i^\g(\overline{\D})\subset U$ for all $\g\in\Ga$ and all $i\in I$,  
as described in Example~\ref{example}. 
For each $\g\in\Ga$ and each $i\in I$, let 
$\tilde{\varphi}_i^\g:=S^{-1}\circ\varphi_i^\g\circ S$ and 
let $\tilde{\Phi}^\g:=\{\tilde{\varphi}_i^\g\}_{i\in I}$. 
Then $\{\tilde{\Phi}^\g\}_{\g\in\Ga}$ is an analytic family 
in $\HIFS(X)$ such that $\g\mapsto\tilde{\Phi}^\g$ is continuous 
with respect to the $\lambda$-topology. 
Since $\{\Phi^\g\}_{\g\in\Ga}$ is of type (V), 
Lemma 3.1 in~\cite{RSU} implies that $\{\tilde{\Phi}^\g\}_{\g\in\Ga}$ is of type (V) as well. 
This completes the proof of our proposition. 
\endpf 
%
%
%

\end{document}